\documentclass[12pt,a4paper,draft]{amsart}

\usepackage{hyperref} 
\usepackage{amssymb}
\usepackage{amsrefs}
\usepackage[all]{xy}

\allowdisplaybreaks

\let\mod=\undefined
\DeclareMathOperator{\Id}{Id} %
\DeclareMathOperator{\Hom}{Hom} %
\DeclareMathOperator{\Ker}{Ker} %
\DeclareMathOperator{\mod}{mod} %
\DeclareMathOperator{\rep}{rep} %
\DeclareMathOperator{\idim}{idim} %
\DeclareMathOperator{\pdim}{pdim} %
\DeclareMathOperator{\proj}{proj} %
\DeclareMathOperator{\umod}{\ul{\mod}} %
\DeclareMathOperator{\urep}{\ul{\rep}} %

\newcommand{\ol}{\overline}
\newcommand{\ul}{\underline}

\newcommand{\bbA}{\mathbb{A}}
\newcommand{\bbN}{\mathbb{N}}
\newcommand{\bbZ}{\mathbb{Z}}

\newcommand{\bone}{\mathbf{1}}
\newcommand{\bGamma}{\mathbf{\Gamma}}

\newcommand{\rt}{\mathrm{t}}
\newcommand{\rs}{\mathrm{s}}

\newcommand{\calA}{\mathcal{A}}
\newcommand{\calD}{\mathcal{D}}
\newcommand{\calI}{\mathcal{I}}
\newcommand{\calJ}{\mathcal{J}}
\newcommand{\calK}{\mathcal{K}}

\newcommand{\definition}{\mathrel{\mathop{=}^{\mathrm{def}}}}
\newcommand{\induction}{\mathrel{\mathop{=}^{\mathrm{ind}}}}

\newcounter{claim}[section]

\newtheorem{corollary}[claim]{Corollary}
\newtheorem{lemma}[claim]{Lemma}
\newtheorem{proposition}[claim]{Proposition}

\newtheorem*{maintheorem}{Main Theorem}
\newtheorem*{fact}{Fact}

\title[Perfect complexes over gentle algebras]{The almost split
triangles for \\ perfect complexes over gentle algebras}

\author{Grzegorz Bobi\'nski}

\address{Faculty of Mathematics and Computer Science \\ Nicolaus
Copernicus University \\ ul.~Chopina 12/18 \\ 87-100 Toru\'n \\
Poland}

\email{gregbob@mat.uni.torun.pl}

\date{}



\begin{document}

\maketitle

Throughout the paper $k$ denotes a fixed field. All vector spaces
and linear maps are $k$-vector spaces and $k$-linear maps,
respectively. By $\bbZ$, $\bbN$, and $\bbN_+$, we denote the sets
of integers, nonnegative integers, and positive integers,
respectively. For $i, j \in \bbZ$, $[i, j] := \{ l \in \bbZ \mid i
\leq l \leq j \}$ (in particular, $[i, j] = \varnothing$ if $i >
j$).

\section*{Introduction}

Given an abelian category $\calA$ one defines its bounded derived
category $\calD^b (\calA)$~\cite{Verdier1977} (of bounded
complexes of objects of $\calA$) having a structure of a
triangulated category, which is an important homological invariant
of $\calA$. In particular, given a finite dimensional algebra $A$
one may study the bounded derived category $\calD^b (\mod A)$ of
the category $\mod A$ of finite dimensional $A$-modules, which one
shortly calls the derived category of $A$ and denotes $\calD^b
(A)$. Since the observation of Happel~\cite{Happel1987}
(generalized by Cline, Parshall, and
Scott~\cite{ClineParshallScott1986}), which states that derived
category is invariant under titling process, an importance of
derived categories in the representation theory of finite
dimensional algebras became clear. This observation was supported
by results connecting derived categories of finite dimensional
algebras with derived categories of coherent sheaves over
projective schemes~\cites{Beilinson1978, GeigleLenzing1987}. Since
that time a lot of results concerning derived categories of finite
dimensional algebras were obtained (see for
example~\cites{Asashiba1999, BocianHolmSkowronski2006,
Brustle2001, Geiss2002, Rickard1989b}). In particular,
Rickard~\cite{Rickard1989} developed the Morita theory for derived
categories of finite dimensional algebras. One of the consequences
is that the derived categories of two finite dimensional algebras
are equivalent as triangulated categories if and only if the
subcategories of perfect complexes are equivalent as triangulated
categories. Recall that if $A$ is a finite dimensional algebra,
then the subcategory of $\calD^b (A)$ formed by perfect complexes
can be identified with the bounded homotopy category $\calK^b
(\proj A)$ of (bounded complexes of) of projective $A$-modules.
Another evidence that the homotopy category of projective modules
delivers a lot of information about the derived category of an
algebra is provided by the result of Krause and
Ye~\cite{KrauseYe2008} stating that the graded centers of $\calK^b
(\proj A)$ and $\calD^b (A)$ are isomorphic for each finite
dimensional algebra $A$.

A class of finite dimensional algebras whose derived categories
attract a lot of interest is the class of gentle algebras
introduced by Assem and Skowro\'nski~\cite{AssemSkowronski1987}.
An important feature of this class of algebra is that it is closed
under derived equivalence, i.e.\ if $A$ is a gentle algebra and
$\calD^b (A)$ is equivalent as a triangulated category to $\calD^b
(B)$ for a finite dimensional algebra $B$, then $B$ is also
gentle~\cite{SchroerZimmermann2003}. Next, this class of algebras
appears naturally in many classification problems. Namely, the
tree gentle algebra are precisely the piecewise hereditary
algebras of type $\bbA$~\cite{AssemHappel1981} (i.e.\ the algebras
derived equivalent to hereditary algebras of type $\bbA$).
Further, if $A$ is a derived discrete algebra, then either $A$ is
piecewise hereditary of Dynkin type or $A$ is a one-cycle gentle
algebra which does not satisfy the clock
condition~\cite{Vossieck2001}. Moreover, the one-cycle gentle
algebras coincide with the piecewise hereditary algebras of type
$\tilde{\bbA}$~\cite{AssemSkowronski1987}. Partial results
concerning the derived equivalence classification of two-cycle
gentle algebras were also obtained~\cites{BobinskiMalicki2008,
AvellaAlaminos2007}. An important role in these investigations is
played by an invariant introduced by Avella-Alaminos and
Geiss~\cite{AvellaAlaminosGeiss2008}.

If $A$ is a gentle algebra, then it is possible to investigate
$\calD^b (A)$ by means of the stable category $\umod \hat{A}$ of
the module category $\mod \hat{A}$ over the repetive algebra
$\hat{A}$~\cite{Ringel1997} (which is no longer finite
dimensional) and the Happel functor~\cite{Happel1988} $\calD^b (A)
\to \umod \hat{A}$. This description is useful, since the
description of the indecomposable objects in $\umod \hat{A}$ is
known. Unfortunately, the precise formula for the Happel functor
seems to be not known. In~\cite{BekkertMerklen2003} Bekkert and
Merklen described the indecomposable objects in $\calD^b (A)$
without using $\hat{A}$, however they did not describe how the
above two descriptions are connected.

Let $A$ be a gentle algebra. Since gentle algebras are
Gorenstein~\cite{GeissReiten2005}, it follows~\cite{Happel1991}
that the almost split triangles in $\calD^b (A)$ exist precisely
for perfect complexes. The aim of this paper is to describe the
almost split triangles in $\calD^b (A)$ in terms of the above
mentioned description of the indecomposable objects in $\calD^b
(A)$ due to Bekkert and Merklen. According to the above remark,
this is equivalent to describing the almost split triangles in
$\calK^b (\proj A)$. The precise formulas are given in
Section~\ref{sect_almostsplit}. The idea of the proof is to use
the Happel functor and the known description of the almost split
triangles in $\umod \hat{A}$. As a side effect we obtain a link
between the two different ways of describing the indecomposable
objects in $\calD^b (A)$.

The paper is organized as follows. In Section~\ref{sect_quivers}
we introduce the language of quivers and their representations,
and in Section~\ref{sect_strings} we present notions of strings
and bands. Next, in Section~\ref{sect_gentle} we present a
description of the indecomposable perfect complexes over gentle
algebras due to Bekker and Merklen, while in
Section~\ref{sect_repetitive} we collect necessary information
about the repetitive algebras of gentle algebras. Finally, in
Section~\ref{sect_Happel} we describe the correspondence between
the indecomposable perfect complexes over a gentle algebra and the
indecomposable modules over its repetitive algebras, and in
Section~\ref{sect_almostsplit} we use this correspondence to
describe the almost split sequences in $\calK^b (\proj A)$.

For basic background on the representation theory of algebras (in
particular, on the tilting theory) we refer
to~\cite{AssemSimsonSkowronski2006}.

The article was written while the author was staying at the
University of Bielefeld as an Alexander von Humboldt Foundation
fellow. The author also acknowledges the support from the Research
Grant No.\ N N201 269135 of the Polish Ministry of Science and
High Education.

\section{Preliminaries on quivers and their representations}
\label{sect_quivers}

By a quiver $\Gamma$ we mean a set $\Gamma_0$ of vertices and a
set $\Gamma_1$ of arrows together with two maps $s, t : \Gamma_1
\to \Gamma_0$ which assign to $\alpha \in \Gamma_1$ the starting
vertex $s \alpha$ and the terminating vertex $t \alpha$,
respectively. We assume that all considered quivers $\Gamma$ are
locally finite, i.e.\ for each $x \in \Gamma_0$ there is only a
finite number of $\alpha \in \Gamma_1$ such that either $s \alpha
= x$ or $t \alpha = x$. A quiver is called finite if $\Gamma_0$
(and, consequently, also $\Gamma_1$) is a finite set. For
technical reasons we assume that all considered quivers $\Gamma$
have no isolated vertices, i.e.\ there is no $x \in \Gamma_0$ such
that $s \alpha \neq x \neq t \alpha$ for each $\alpha \in
\Gamma_1$.

Let $\Gamma$ be a quiver. If $l \in \bbN_+$, then by a path in
$\Gamma$ of length $l$ we mean $\sigma = \alpha_1 \cdots \alpha_l$
such that $\alpha_i \in \Gamma_1$ for each $i \in [1, l]$ and $s
\alpha_i = t \alpha_{i + 1}$ for each $i \in [1, l - 1]$. In the
above situation we put $s \sigma := s \alpha_l$ and $t \sigma := t
\alpha_1$. Moreover, we put $\alpha_i (\sigma) := \alpha_i$ for $i
\in [1, l]$. Observe that each $\alpha \in \Gamma$ is a path in
$\Gamma$ of length $1$. Moreover, for each $x \in \Gamma_0$ we
introduce the path $\bone_x$ in $\Gamma$ of length $0$ such that
$s \bone_x := x =: t \bone_x$. We denote the length of a path
$\sigma$ in $\Gamma$ by $\ell (\sigma)$. If $\sigma'$ and
$\sigma''$ are two paths in $\Gamma$ such that $s \sigma' = t
\sigma''$, then we define the composition $\sigma' \sigma''$ of
$\sigma'$ and $\sigma''$, which is a path in $\Gamma$ of length
$\ell (\sigma') + \ell (\sigma'')$, in the obvious way (in
particular, $\sigma \bone_{s \sigma} = \sigma = \bone_{t \sigma}
\sigma$ for each path $\sigma$). In order to increase clarity we
sometimes write $\sigma' \cdot \sigma''$ instead of $\sigma'
\sigma''$ in the above situation. If $\sigma$ is a path such that
$s \sigma = t \sigma$, then for $n \in \bbN_+$ we denote by
$\sigma^n$ the $n$-fold composition of $\sigma$ with itself.

Let $\Gamma$ be a quiver. We define the double quiver
$\ol{\Gamma}$ of $\Gamma$ in the following way: $\ol{\Gamma}_0 :=
\Gamma_0$, $\ol{\Gamma}_1 := \Gamma_1 \cup \Gamma_1^{-1}$, where
$\Gamma_1^{-1} := \{ \alpha^{-1} \mid \alpha \in \Gamma_1 \}$, and
$s \alpha^{-1} := t \alpha$ and $t \alpha^{-1} := s \alpha$ for
$\alpha \in \Gamma_1$. By $\approx$ we denote the equivalence
relation in $\ol{\Gamma}_1$ whose residue classes are $\Gamma_1$
and $\Gamma_1^{-1}$. We put $(\alpha^{-1})^{-1} := \alpha$ for
$\alpha \in \Gamma_1$ and extend the operation $(-)^{-1}$ to the
paths in $\ol{\Gamma}$ of positive length in such a way that
$(\omega' \omega'')^{-1} = \omega''^{-1} \omega'^{-1}$ for all
paths $\omega'$ and $\omega''$ in $\ol{\Gamma}$ of positive length
such that $s \omega' = t \omega''$. If $\omega$ is a path in
$\ol{\Gamma}$ of positive length and $i \in [1, \ell (\omega)]$,
then $\alpha_i^{-1} (\omega) := (\alpha_i (\omega))^{-1}$. For a
set $\Sigma$ of paths in $\ol{\Gamma}$ of positive length we put
$\Sigma^{-1} := \{ \sigma^{-1} \mid \sigma \in \Sigma \}$.

Let $\Gamma$ be a quiver. We define the path category $k \Gamma$
of $\Gamma$ as follows. The objects of $k \Gamma$ are the vertices
of $\Gamma$. If $x', x'' \in \Gamma_0$, then the homomorphism
space $k \Gamma (x', x'')$ consists of the formal $k$-linear
combinations of paths starting at $x'$ and terminating at $x''$.
The composition of maps in $k \Gamma$ is induced by the
composition of paths in $\Gamma$. For a set $R$ of morphisms in $k
\Gamma$ we denote by $\langle R \rangle$ the ideal in $k \Gamma$
generated by $R$. A morphism $\varrho$ in $\Gamma$ is called a
relation if $\varrho \in \langle \Gamma_1 \rangle^2$. A set $R$ of
relations in $k \Gamma$ is called admissible if there exists $n
\in \bbN_+$ such that $\langle \Gamma_1 \rangle^n \subset \langle
R \rangle$.

By an (admissible) bound quiver we mean a pair $(\Gamma, R)$
consisting of a quiver $\Gamma$ and an (admissible, respectively)
set of relations in $k \Gamma$. For a bound quiver $\bGamma =
(\Gamma, R)$ we denote by $k \bGamma$ the corresponding factor
category $k \Gamma / \langle R \rangle$. If $\bGamma = (\Gamma,
R)$ is a bound quiver and $\varrho \in k \bGamma (x', x'')$ for
$x', x'' \in \Gamma_0$, then we put $s \varrho := x'$ and $t
\varrho := x''$. A bound quiver $(\Gamma, R)$ is called monomial
if $R$ consists of paths.

Let $\bGamma = (\Gamma, R)$ be a monomial bound quiver. By a path
in $\bGamma$ we mean a path in $\Gamma$ which does not belong to
$\langle R \rangle$. If $x', x'' \in \Gamma_0$, then we identify
$k \bGamma (x', x'')$ with the subspace of $k \Gamma (x', x'')$
spanned by the paths in $\bGamma$ starting at $x'$ and terminating
at $x''$. A path $\sigma$ in $\bGamma$ is said to be maximal in
$\bGamma$ if there are no paths $\sigma'$ and $\sigma''$ in
$\bGamma$ such that $s \sigma' = t \sigma$, $t \sigma'' = s
\sigma'$, $\sigma' \sigma \sigma''$ is a path in $\bGamma$, and
$\ell (\sigma') + \ell (\sigma'') > 0$. The lack of isolated
vertices in $\Gamma$ implies that $\ell (\sigma) > 0$ for each
maximal path $\sigma$ in $\bGamma$.

For the rest of the section we assume that $\bGamma = (\Gamma, R)$
is an admissible bound quiver.

By a representation of $\bGamma$ we mean a functor $M : k \bGamma
\to \mod k$, where $\mod k$ denotes the category of finite
dimensional vector spaces over $k$, such that $M (x) \neq 0$ only
for a finite number of $x \in \Gamma_0$. Observe that a
representation $M$ of $\bGamma$ is uniquely determined by the
collection $(M (x))_{x \in \Gamma_0}$ of vector spaces and the
collection $(M (\alpha))_{\alpha \in \Gamma_1}$ of linear maps. On
the other hand, a pair of such collections determines a
representation of $\bGamma$ if and only if the induced map $M
(\varrho)$ vanishes for each $\varrho \in R$. If $M$ and $N$ are
two representations of $\bGamma$, then the morphism space
$\Hom_{\bGamma} (M, N)$ consists of the natural transformations of
the corresponding functors. We denote the category of
representations of $\bGamma$ by $\rep \bGamma$. It is well-known
that $\rep \bGamma$ is an abelian category which possesses almost
split sequences. We denote by $\tau_{\bGamma}$ the
Auslander--Reiten translation in $\rep \bGamma$. We remark that
Gabriel proved (see for
example~\cite{AssemSimsonSkowronski2006}*{Corollaries~I.6.10
and~II.3.7}) that for each finite dimensional algebra $A$ the
category of $A$-modules is equivalent to the category of
representations for an appropriate admissible bound quiver. This
implies in particular, that we may work with bound quivers instead
of algebras.

Now we describe the indecomposable projective representations of
$\bGamma$. For each $x \in \Gamma_0$ we define $P_x \in \rep
\bGamma$ as follows: $P_x (x') := k \bGamma (x, x')$ for $x' \in
\Gamma_0$ and $P_x (\varrho) (\varrho') := \varrho \varrho'$ for
morphisms $\varrho$ and $\varrho'$ in $k \bGamma$ such that $s
\varrho' = x$ and $s \varrho = t \varrho'$. Moreover, if $\varrho$
is a morphism in $k \bGamma$, then we define $p_\varrho : P_{t
\varrho} \to P_{s \varrho}$ by $p_\varrho (x) (\varrho') :=
\varrho' \varrho$ for $x \in \Gamma_0$ and a morphism $\varrho'$
in $k \bGamma$ with $s \varrho' = t \varrho$ and $t \varrho' = x$.
It is an easy exercise to check that the map $\Hom_{\bGamma} (P_x,
M) \to M (x)$, $f \mapsto f (x) (\bone_x)$, is an isomorphism of
vector spaces for each $x \in \Gamma_0$ and $M \in \rep \bGamma$.
This implies that the above formulas describe the fully faithful
contravariant functor $\bGamma \to \rep \bGamma$ whose essential
image coincides with the full subcategory of indecomposable
projective representations of $\bGamma$.

Similarly, we describe the indecomposable injective
representations of $\bGamma$. For $x \in \Gamma_0$ we define $Q_x
\in \rep \bGamma$ by $Q_x (x') := (k \bGamma (x', x))^*$ for $x'
\in \Gamma_0$, where $(-)^* : \mod k \to \mod k$ denotes the
$k$-linear dual, and $Q_x (\varrho) (\varphi) (\varrho') :=
\varphi (\varrho' \varrho)$ for morphisms $\varrho$ and $\varrho'$
in $k \bGamma$ such that $s \varrho' = t \varrho$ and $t \varrho'
= x$, and $\varphi \in (k \bGamma (s \varrho, x))^*$. Moreover, if
$\varrho$ is a morphism in $k \bGamma$, then we define $q_\varrho
: Q_{t \varrho} \to Q_{s \varrho}$ by $q_\varrho (x) (\varphi)
(\varrho') := \varphi (\varrho \varrho')$ for $x \in \Gamma_0$,
morphisms $\varrho$ and $\varrho'$ in $k \bGamma$ such that $s
\varrho' = x$ and $t \varrho' = s \varrho$, and $\varphi \in (k
\bGamma (x, t \varrho))^*$. Again, the map $(M (x))^* \to
\Hom_{\bGamma} (M, Q_x)$, $\varphi \mapsto (m \mapsto (\varrho \to
\varphi (M (\varrho) (m))))$, is an isomorphism for each $x \in
\Gamma_0$ and $M \in \rep \bGamma$, and, consequently, we obtain
the fully faithful contravariant functor $\bGamma \to \rep
\bGamma$ whose essential image coincides with the full subcategory
of indecomposable injective representations of $\bGamma$.

\section{Almost gentle quivers} \label{sect_strings}

An admissible monomial bound quiver $\bGamma = (\Gamma, R)$ is
called almost gentle if the following conditions are satisfied:
\begin{enumerate}

\item
for each $x \in \Gamma_0$ there are at most two $\alpha \in
\Gamma_1$ such that $s \alpha = x$ and at most two $\alpha \in
\Gamma_1$ such that $t \alpha = x$,

\item
for each $\alpha \in \Gamma_1$ there is at most one $\alpha' \in
\Gamma_1$ such that $s \alpha' = t \alpha$ and $\alpha' \alpha
\not \in R$, and at most one $\alpha' \in \Gamma_1$ such that $t
\alpha' = s \alpha$ and $\alpha \alpha' \not \in R$,

\item
for each $\alpha \in \Gamma_1$ there is at most one $\alpha' \in
\Gamma_1$ such that $s \alpha' = t \alpha$ and $\alpha' \alpha \in
R$, and at most one $\alpha' \in \Gamma_1$ such $t \alpha' = s
\alpha$ and $\alpha \alpha' \in R$.

\end{enumerate}
Equivalently, $\bGamma$ is an almost gentle quiver if and only if
there exist functions $S, T : \Gamma_1 \to \{ \pm 1 \}$, which we
call string functions for $\bGamma$, such that the following
conditions are satisfied:
\begin{enumerate}

\item
if $s \alpha' = s \alpha''$ and $\alpha' \neq \alpha''$ for
$\alpha', \alpha'' \in \Gamma_1$, then $S \alpha' = - S \alpha''$,

\item
if $t \alpha' = t \alpha''$ and $\alpha' \neq \alpha''$ for
$\alpha', \alpha'' \in \Gamma_1$, then $T \alpha' = - T \alpha''$,

\item
if $s \alpha' = t \alpha''$ and $\alpha' \alpha'' \not \in R$ for
$\alpha', \alpha'' \in \Gamma_1$, then $S \alpha' = - T \alpha''$.

\item
if $s \alpha' = t \alpha''$ and $\alpha' \alpha'' \in R$ for
$\alpha', \alpha'' \in \Gamma_1$, then $S \alpha' = T \alpha''$.

\end{enumerate}
Note that the string functions for $\bGamma$ are not uniquely
determined by $\bGamma$. For the rest of the section we fix an
almost gentle bound quiver $\bGamma = (\Gamma, R)$ together with
string functions $S$ and $T$.

Let $R' := R \cup R^{-1} \cup \{ \alpha \alpha^{-1}, \alpha^{-1}
\alpha \mid \alpha \in \Gamma_1 \}$. Then $(\ol{\Gamma}, R')$ is a
monomial bound quiver. If $l \in \bbN_+$, then by a string in
$\bGamma$ of lenght $l$ we mean a path in $(\ol{\Gamma}, R')$ of
length $l$. Moreover, for each $x \in \Gamma_0$ we introduce two
strings $\bone_{x, 1}$ and $\bone_{x, -1}$ such that $\ell
(\bone_{x, \varepsilon}) := 0$ and $s \bone_{x, \varepsilon} := x
=: t \bone_{x, \varepsilon}$ for $\varepsilon \in \{ \pm 1 \}$. We
put $(\bone_{x, \varepsilon})^{-1} := \bone_{x, -\varepsilon}$ for
$x \in \Gamma_0$ and $\varepsilon \in \{ \pm 1 \}$. Observe that
every path in $\bGamma$ of positive length is a string in
$\bGamma$. A string $\omega$ in $\bGamma$ is called simple if
either $\ell (\omega) = 0$ or $\omega$ is a path in $\bGamma$ (of
positive length). Moreover, a string $\omega$ in $\bGamma$ is
called directed if either $\omega$ or $\omega^{-1}$ is a simple
string. A string $\omega$ in $\bGamma$ is called a band if $\ell
(\omega) > 0$, either $\alpha_1 (\omega) \in \Gamma_1$ and
$\alpha_{\ell (\omega)}^{-1} (\omega) \in \Gamma_1$ or
$\alpha_1^{-1} (\omega) \in \Gamma_1$ and $\alpha_{\ell (\omega)}
(\omega) \in \Gamma_1$, $s \omega = t \omega$, $\omega^n$ is a
string in $\bGamma$ for each $n \in \bbN_+$, and there is no
string $\omega'$ in $\bGamma$ such that $\ell (\omega') < \ell
(\omega)$, $s \omega' = t \omega'$, and $\omega = \omega'^n$ for
some $n \in \bbN_+$.

We extend the functions $S$ and $T$ to the strings in $\bGamma$ as
follows. First, we put $S \alpha^{-1} := T \alpha$ and $T
\alpha^{-1} := S \alpha$ for $\alpha \in \Gamma_1$. Next, we put
$S \omega := S \alpha_{\ell (\omega)} (\omega)$ and $T \omega := T
\alpha_1 (\omega)$ for a string $\omega$ in $\bGamma$ of positive
length. Finally, we put $S \bone_{x, \varepsilon} := \varepsilon$
and $T \bone_{x, \varepsilon} := - \varepsilon$ for $x \in
\Gamma_0$ and $\varepsilon \in \{ \pm 1 \}$. Observe that if
$\omega'$ and $\omega''$ are strings in $\bGamma$ of positive
length such that $s \omega' = t \omega''$ and $\omega' \omega''$
is a string in $\bGamma$, then $S \omega' = - T \omega''$.
Consequently, if $\omega$ is a string in $\bGamma$, $x \in
\Gamma_0$, and $\varepsilon \in \{ \pm 1 \}$, then we say that the
composition $\omega \bone_{x, \varepsilon}$ ($\bone_{x,
\varepsilon} \omega$) is defined (and equals $\omega$) if and only
if $x = s \omega$ and $\varepsilon = S \omega$ ($x = t \omega$ and
$\varepsilon = - T \omega$, respectively).

Let $\omega$ be a string in $\bGamma$ of length $l$. If $i \in [0,
l]$, then we denote by $\omega_{[i]}$ and ${}_{[i]} \omega$ the
strings in $\bGamma$ of length $i$ and $l - i$, respectively, such
that $\omega = \omega_{[i]} \cdot {}_{[i]} \omega$. In particular,
$\omega_{[0]} = \bone_{t \omega, - T \omega}$ and $_{[\ell
(\omega)]} \omega = \bone_{s \omega, S \omega}$.

Fix $x \in \Gamma_0$ and $\varepsilon \in \{ \pm 1 \}$. By
$\Sigma_{x, \varepsilon}$ we denote the set of simple strings
$\sigma$ in $\bGamma$ such that $s \sigma = x$ and $S \sigma =
\varepsilon$. Similarly, $\Sigma_{x, \varepsilon}'$ denotes the
set of simple strings $\sigma$ in $\bGamma$ such that $t \sigma =
x$ and $T \sigma = \varepsilon$. Next,
\[
\alpha_{x, \varepsilon} :=
\begin{cases}
\alpha & \alpha \in \Sigma_{x, \varepsilon} \cap \Gamma_1,
\\
\varnothing & \text{otherwise},
\end{cases}
\qquad \text{and} \qquad \alpha_{x, \varepsilon}' :=
\begin{cases}
\alpha & \alpha \in \Sigma_{x, \varepsilon}' \cap \Gamma_1,
\\
\varnothing & \text{otherwise}.
\end{cases}
\]
Finally, we denote by $\sigma_{x, \varepsilon}$ and $\sigma_{x,
\varepsilon}'$ the strings of maximal length in $\Sigma_{x,
\varepsilon}$ and $\Sigma_{x, \varepsilon}'$, respectively.

\section{The homotopy category of a gentle quiver}
\label{sect_gentle}

An almost gentle bound quiver $\bGamma = (\Gamma, R)$ is called
gentle if $\Gamma$ is finite and $R$ consists of paths of length
$2$. For the rest of the section we assume that $\bGamma =
(\Gamma, R)$ is a fixed gentle bound quiver together with string
functions $S$ and $T$.

Let $\calK^b (\bGamma)$ denote the bounded homotopy category of
complexes of projective representations of $\bGamma$. Recall that
$\calK^b (\bGamma)$ has a structure of a triangulated
category~\cite{KonigZimmermann1998}*{Theorem~2.3.1} with the
suspension functor given by the degree shift $X \mapsto X [1]$.
Moreover, since $\bGamma$ is Gorenstein, i.e.\ $\pdim_{\bGamma} Q
< \infty$ for each injective representation $Q$ of $\bGamma$ and
$\idim_{\bGamma} P < \infty$ for each projective representation
$P$ of $\bGamma$~\cite{GeissReiten2005}, $\calK^b (\bGamma)$
possesses almost split
triangles~\cite{HappelKellerReiten2008}*{Section~5}, thus also the
Auslander--Reiten translation, which we denote by $\tau_{\calK^b
(\bGamma)}$.

Let $R' := \{ \alpha \alpha^{-1}, \alpha^{-1} \alpha \mid \alpha
\in \Gamma_1 \}$. Then $(\ol{\Gamma}, R')$ is a monomial bound
quiver. If $l \in \bbN_+$, then by a homotopy string in $\bGamma$
of length $l$ we mean a path in $(\ol{\Gamma}, R')$ of length $l$.
Moreover, $\bone_{x, 1}$ and $\bone_{x, -1}$ are homotopy strings
in $\bGamma$ of length $0$ for each $x \in \Gamma_0$. Observe that
every string in $\bGamma$ is a homotopy string in $\bGamma$. We
extend $S$ and $T$ to the homotopy strings in $\bGamma$ in the
usual way. If $\omega'$ and $\omega''$ are homotopy strings in
$\bGamma$ of positive length, then we say that the composition
$\omega' \omega''$ is defined (in the obvious way) if $s \omega' =
t \omega''$ and one of the following conditions is satisfied,
where $\alpha' := \alpha_{\ell (\omega')} (\omega')$ and $\alpha''
:= \alpha_1 (\omega'')$:
\begin{enumerate}

\item
$S \omega' = T \omega''$ and either $\alpha', \alpha'' \in
\Gamma_1$ or $\alpha'^{-1}, \alpha''^{-1} \in \Gamma_1$,

\item
$S \omega' = - T \omega''$ and either $\alpha', \alpha''^{-1} \in
\Gamma_1$ or $\alpha'^{-1}, \alpha'' \in \Gamma_1$.

\end{enumerate}
Similarly, if $\omega$ is a homotopy string in $\bGamma$ of
positive length, $x \in \Gamma_0$, and $\varepsilon \in \{ \pm 1
\}$, then the composition $\omega \bone_{x, \varepsilon}$
($\bone_{x, \varepsilon} \omega$) is defined (and equals $\omega$)
if and only if $x = s \omega$ and either $\varepsilon = S \omega$
and $\alpha_{\ell (\omega)} (\omega) \in \Gamma_1$ or $\varepsilon
= - S \omega$ and $\alpha_{\ell (\omega)}^{-1} (\omega) \in
\Gamma_1$ ($x = t \omega$ and either $\varepsilon = T \omega$ and
$\alpha_1 (\omega) \in \Gamma_1$ or $\varepsilon = - T \omega$ and
$\alpha_1^{-1} (\omega) \in \Gamma_1$, respectively). Finally, if
$x', x'' \in \Gamma_0$ and $\varepsilon', \varepsilon'' \in \{ \pm
1 \}$, then the composition $\bone_{x', \varepsilon'} \bone_{x'',
\varepsilon''}$ is defined (and equals $\bone_{x', \varepsilon'}$)
if and only if $x' = x''$ and $\varepsilon' = \varepsilon''$.
Observe that the above definitions for homotopy strings differ
from the ones we have for strings. If $\omega$ is a homotopy
string in $\bGamma$, then by $\sigma_\omega$ we denote the string
of maximal length among the simple strings $\sigma$ in $\bGamma$
such that the composition $\sigma \omega$ (as homotopy strings in
$\bGamma$) is defined.

A simple homotopy string $\theta$ in $\bGamma$ is called an
antipath in $\bGamma$ provided $\alpha_i (\theta) \alpha_{i + 1}
(\theta) \in R$ (equivalently, $S \alpha_i (\theta) = T \alpha_{i
+ 1} (\theta)$) for each $i \in [1, \ell (\theta) - 1]$. For $x
\in \Gamma_0$ and $\varepsilon \in \{ \pm 1 \}$, let $\Theta_{x,
\varepsilon}$ denote the set of all antipaths $\theta$ in
$\bGamma$ such that $t \theta = x$ and $T \theta = \varepsilon$.
If there is an antipath in $\Theta_{x, \varepsilon}$ of maximal
length, then we denote it by $\theta_{x, \varepsilon}$. Otherwise,
we put $\theta_{x, \varepsilon} := \varnothing$.

Let $\omega$ be a homotopy string in $\bGamma$. If $\ell (\omega)
> 0$, then $\omega$ has a unique presentation in the form $\omega
= \sigma_1 \cdots \sigma_L$, $L \in \bbN_+$, such that $\sigma_i$
is a directed string in $\bGamma$ of positive length for each $i
\in [1, L]$, and the composition of $\sigma_i \sigma_{i + 1}$ (as
homotopy strings in $\bGamma$) is defined for each $i \in [1, L -
1]$. In the above situation we put $L (\omega) := L$, $\sigma_i
(\omega) := \sigma_i$ and $\sigma_i^{-1} (\omega) :=
\sigma_i^{-1}$ for $i \in [1, L]$, and
\begin{multline*}
\deg \omega := |\{ i \in [1, L] \mid \text{$\sigma_i$ is a path in
$\bGamma$} \}|
\\
- |\{ i \in [1, L] \mid \text{$\sigma_i^{-1}$ is a path in
$\bGamma$} \}|.
\end{multline*}
Moreover, we put $L (\omega) := 0$ and $\deg \omega := 0$ if $\ell
(\omega) = 0$. If $i \in [0, L (\omega)]$, then we denote by
$\omega^{[i]}$ and ${}^{[i]} \omega$ the homotopy strings in
$\bGamma$ of length $\sum_{j \in [1, i]} \ell (\sigma_j (\omega))$
and $\sum_{j \in [i +  1, L (\omega)]} \ell (\sigma_j (\omega))$,
respectively, such that $\omega = \omega^{[i]} \cdot {}^{[i]}
\omega$. In particular,
\begin{align*}
\omega^{[0]} & =
\begin{cases}
\bone_{t \omega, T \omega} & \text{$\ell (\omega) > 0$ and
$\alpha_1 (\omega) \in \Gamma_1$},
\\
\bone_{t \omega, - T \omega} & \text{otherwise},
\end{cases}
\\
\intertext{and} %
{}^{[L (\omega)]} \omega & =
\begin{cases}
\bone_{s \omega, - S \omega} & \text{$\ell (\omega) > 0$ and
$\alpha_{\ell (\omega)}^{-1} (\omega) \in \Gamma_1$},
\\
\bone_{t \omega, S \omega} & \text{otherwise}.
\end{cases}
\end{align*}
Moreover, ${}^{[i]} \omega^{[j]} := {}^{[i]} (\omega^{[j]})$ for
$i, j \in [0, L (\omega)]$, $i \leq j$.

Let $\omega$ be a homotopy string in $\bGamma$ and $m \in \bbZ$.
We define $X = X_{m, \omega} \in \calK^b (\bGamma)$ in the
following way. First, for $m' \in \bbZ$ we put $\calI_{m'} =
\calI_{m'} (m, \omega) := \{ i \in [0, L (\omega)] \mid m + \deg
\omega^{[i]} = m' \}$. Then $X^{m'} := \bigoplus_{i \in
\calI_{m'}} P_{s \omega^{[i]}}$ for $m' \in \bbZ$ and
\[
(d_X^{m'})_{i, j} :=
\begin{cases}
p_{\sigma_{j + 1} (\omega)} & \text{$i = j + 1$ and $\sigma_{j +
1} (\omega)$ is a path in $\bGamma$},
\\
p_{\sigma_j^{-1} (\omega)} & \text{$i = j - 1$ and $\sigma_j^{-1}
(\omega)$ is a path in $\bGamma$},
\\
0 & \text{otherwise},
\end{cases}
\]
for $m' \in \bbZ$, $j \in \calI_{m'}$, and $i \in \calI_{m' + 1}$.
The objects of $\calK^b (\bGamma)$ of the above form are called
the string complexes.

For a homotopy string $\omega$ in $\bGamma$ and $m \in \bbZ$ we
denote by $\Upsilon_{m, \omega}$ the map $\Upsilon : X_{m, \omega}
\to X_{m + \deg \omega, \omega^{-1}}$ defined by
\[
\Upsilon_{i, j}^{m'} :=
\begin{cases}
p_{s \omega^{[j]}} & \text{$i = L (\omega) - j$},
\\
0 & \text{otherwise},
\end{cases}
\]
for $m' \in \bbZ$, $j \in \calI_{m'} (m, \omega)$, and $i \in
\calI_{m'} (m + \deg \omega, \omega^{-1})$. Observe that
$\Upsilon_{m, \omega}$ is an isomorphism for each homotopy string
$\omega$ in $\bGamma$ and $m \in \bbZ$ --- the inverse map is
given by $\Upsilon_{m + \deg \omega, \omega^{-1}}$.

Let $\omega$ be a homotopy string in $\bGamma$. Let $\sigma$ be a
path in $\bGamma$ of positive length such that the composition
$\sigma \omega$ (as homotopy strings in $\bGamma$) is defined. For
$m \in \bbZ$ we denote by $F_{m, \sigma, \omega}'$ the map $F' :
X_{m, \bone_{t \sigma, - T \sigma}} \to X_{m, \omega}$ defined by
\[
F_{i, 0}'^m :=
\begin{cases}
p_\sigma & \text{$i = 0$},
\\
0 & \text{otherwise},
\end{cases}
\]
$i \in \calI_m (m, \omega)$. Similarly, let $\sigma$ be a path in
$\bGamma$ of positive length such that the composition
$\sigma^{-1} \omega$ is defined. We denote by $F_{m, \sigma,
\omega}''$ the map $F'' : X_{m, \omega} \to X_{m, \bone_{s \sigma,
S \sigma}}$ defined by
\[
F_{0, j}'^m :=
\begin{cases}
p_\sigma & \text{$j = 0$},
\\
0 & \text{otherwise},
\end{cases}
\]
for $j \in \calI_m (m, \omega)$. Next, if $m \in \bbZ$ and
$\sigma$ is a path in $\bGamma$ of positive length such that the
composition $\omega \sigma^{-1}$ is defined, then we put $G_{m,
\sigma, \omega}' := \Upsilon_{m + \deg \omega, \omega^{-1}} \circ
F_{m + \deg \omega, \sigma, \omega^{-1}}' \circ \Upsilon_{m + \deg
\omega, \bone_{t \sigma, T \sigma}} : X_{m + \deg \omega, \bone_{t
\sigma, T \sigma}} \to X_{m, \omega}$. Finally, if $m \in \bbZ$
and $\sigma$ is a path in $\bGamma$ of positive length such that
the composition $\omega \sigma$ is defined, then we put $G_{m,
\sigma, \omega}'' := \Upsilon_{m + \deg \omega, \bone_{s \sigma, S
\sigma}} \circ F_{m + \deg \omega, \sigma, \omega^{-1}}'' \circ
\Upsilon_{m, \omega} : X_{m, \omega} \to X_{m + \deg \omega,
\bone_{s \sigma, - S \sigma}}$.

A homotopy string $\omega$ in $\bGamma$ is called a homotopy band
if $\deg \omega = 0$, $L (\omega) > 0$, either $\sigma_1 (\omega)$
and $\sigma_{L (\omega)}^{-1} (\omega)$ are paths in $\bGamma$ or
$\sigma_1^{-1} (\omega)$ and $\sigma_{L (\omega)} (\omega)$ are
paths in $\bGamma$, $s \omega = t \omega$, $S \omega = - T
\omega$, and there is no homotopy string $\omega'$ in $\bGamma$
such that $\ell (\omega') < \ell (\omega)$, $s \omega' = t
\omega'$, and $\omega = \omega'^n$ for some $n \in \bbN_+$. If
$\omega$ is a homotopy band in $\bGamma$ and $i \in [0, L (\omega)
- 1]$, then we put $\omega^{(i)} := {}^{[i]} \omega \cdot
\omega^{[i]}$. Observe that it may happen that $\omega^{(i)}$ is
not a homotopy band in the above situation.

Let $\omega$ be a homotopy band in $\bGamma$, $\mu$ an
indecomposable automorphism of a finite dimensional vector space
$K$, and $m \in \bbZ$. We define $Y = Y_{m, \omega, \mu} \in
\calK^b (\bGamma)$ in the following way. First, for $m' \in \bbZ$
we put $\calJ_{m'} := \{ i \in [1, L (\omega)] \mid m + \deg
\omega^{[i]} = m' \}$. Then $Y^{m'} := \bigoplus_{i \in
\calJ_{m'}} P_{s \omega^{[i]}} \otimes_k K$ for $m' \in \bbZ$ and
\[
(d_Y^{m'})_{i, j} :=
\begin{cases}
p_{\sigma_i (\omega)} \otimes \Id & \text{$i = j + 1$, $i > 1$,}
\\
& \quad \text{and $\sigma_{j + 1} (\omega)$ is a path in
$\bGamma$},
\\
p_{\sigma_j^{-1} (\omega)} \otimes \Id & \text{$i = j - 1$, $j < L
(\omega)$},
\\
& \quad \text{and $\sigma_j^{-1} (\omega)$ is a path in
$\bGamma$},
\\
p_{\sigma_1 (\omega)} \otimes \mu & \text{$j = L (\omega)$, $i =
1$, $L (\omega) > 2$},
\\
& \quad \text{and $\sigma_1 (\omega)$ is a path in $\bGamma$},
\\
p_{\sigma_1^{-1} (\omega)} \otimes \mu & \text{$j = 1$, $i = L
(\omega)$, $L (\omega) > 2$},
\\
& \quad \text{and $\sigma_1^{-1} (\omega)$ is a path in
$\bGamma$},
\\
p_{\sigma_1 (\omega)} \otimes \mu + p_{\sigma_2^{-1} (\omega)}
\otimes \Id & \text{$j = L (\omega)$, $i = 1$, $L (\omega) = 2$},
\\
& \quad \text{and $\sigma_1 (\omega)$ is a path in $\bGamma$},
\\
p_{\sigma_1^{-1} (\omega)} \otimes \mu + p_{\sigma_2 (\omega)}
\otimes \Id & \text{$j = 1$, $i = L (\omega)$, $L (\omega) = 2$},
\\
& \quad \text{and $\sigma_1^{-1} (\omega)$ is a path in
$\bGamma$},
\\
0 & \text{otherwise},
\end{cases}
\]
for $m' \in \bbZ$, $j \in \calJ_{m'}$, and $i \in \calJ_{m' + 1}$.
The objects of $\calK^b (\bGamma)$ of the above form are called
the band complexes.

The following description of the indecomposable objects in
$\calK^b (\bGamma)$ was obtained in~\cite{BekkertMerklen2003}.

\begin{proposition}
Let $\bGamma$ be a gentle bound quiver. If $X$ is an
indecomposable object in $\calK^b (\bGamma)$, then either $X
\simeq X_{m, \omega}$ for some $m \in \bbZ$ and a homotopy string
$\omega$ in $\bGamma$ or $X \simeq Y_{m, \omega, \mu}$ for some $m
\in \bbZ$, a homotopy band $\omega$ in $\bGamma$, and an
automorphism $\mu$ of a finite dimensional vector space.
\end{proposition}

\section{The repetitive quiver of a gentle quiver}
\label{sect_repetitive}

Throughout this section we fix a gentle bound quiver $\bGamma  =
(\Gamma, R)$ together with string functions $S$ and $T$. We also
denote by $\Sigma$ the set of maximal paths in $\bGamma$.

Our first aim in this section is to define the repetitive quiver
$\hat{\bGamma} = (\hat{\Gamma}, \hat{R})$ of $\bGamma$. We put
$\hat{\Gamma}_0 := \Gamma_0 \times \bbZ$ and $\hat{\Gamma}_1 :=
(\Gamma_1 \times \bbZ) \cup (\Sigma^* \times \bbZ)$, where
$\Sigma^* := \{ \sigma^* \mid \sigma \in \Sigma \}$. Moreover, $s
(\alpha [m]) := (s \alpha) [m]$ and $t (\alpha [m]) := (t \alpha)
[m]$ for $\alpha \in \Gamma_1$ and $m \in \bbZ$, and $s (\sigma^*
[m]) := (t \sigma) [m + 1]$ and $t (\sigma^* [m]) := (s \sigma)
[m]$ for $\sigma \in \Sigma$ and $m \in \bbZ$. For a path $\sigma$
in $\Gamma$ and $m \in \bbZ$ we define the path $\sigma [m]$ in
$\hat{\Gamma}$ in the obvious way. Let
\begin{align*}
\bbZ R := & \{ \sigma [m] \mid \sigma \in R, \; m \in \bbZ \} \cup
\\
& \{ \sigma'^* [m - 1] \sigma''^* [m] \mid \sigma', \sigma'' \in
\Sigma, \; t \sigma' = s \sigma'', \; m \in \bbZ \} \cup
\\
& \{ \alpha [m] \sigma^* [m] \mid \alpha \in \Gamma_1, \; \sigma
\in \Sigma, \; s \alpha = s \sigma, \; S \alpha = - S \sigma , \;
m \in \bbZ \} \cup
\\
& \{ \sigma^* [m] \alpha [m + 1] \mid \alpha \in \Gamma_1, \;
\sigma \in \Sigma, \; t \alpha = t \sigma, \; T \alpha = - T
\sigma, \; m \in \bbZ \}.
\end{align*}
Every path in $\hat{\Gamma}$ of the form $\sigma'' [m] \sigma^*
[m] \sigma' [m + 1]$, where $\sigma \in \Sigma$, $\sigma = \sigma'
\sigma''$ for paths $\sigma'$ and $\sigma''$ in $\Gamma$, and $m
\in \bbZ$, is called a full path. Let $\Lambda$ be the set of full
paths in $\hat{\Gamma}$. Then
\begin{align*}
\hat{R} := & \bbZ R \cup \{ \lambda' - \lambda'' \mid \lambda',
\lambda'' \in \Lambda, \; s \lambda' = s \lambda'', \; \lambda'
\neq \lambda'' \} \cup
\\
& \{ \beta \lambda \mid \beta \in \hat{\Gamma}_1, \; \lambda \in
\Lambda, \; s \beta = t \lambda \} \cup \{ \lambda \beta \mid
\beta \in \hat{\Gamma}_1, \; \lambda \in \Lambda, \; s \lambda = t
\beta \}.
\end{align*}
We have the Nakayama automorphism $\nu$ of $\hat{\bGamma}$ given
by $\nu (x [m]) := x [m + 1]$ and $\nu (\alpha [m]) := \alpha [m +
1]$ for $x \in \Gamma_0$, $\alpha \in \Gamma_1 \cup \Sigma^*$, and
$m \in \bbZ$.

Let $\urep \hat{\bGamma}$ denote the stable category of $\bGamma$.
Since $\rep \hat{\bGamma}$ is a Frobenius category, $\urep
\hat{\bGamma}$ is a triangulated category with the suspension
functor given by the inverse $\Omega^{-1}$ of the Heller syzygy
functor $\Omega$. Every exact sequence in $\rep \hat{\bGamma}$
induces a triangle in $\urep \hat{\bGamma}$. Moreover, $\urep
\hat{\bGamma}$ possesses almost split triangles, which come from
the almost split sequences in $\rep \hat{\bGamma}$. In particular,
the Auslander--Reiten translation in $\urep \hat{\bGamma}$ is
given by the Auslander--Reiten translation $\tau_{\hat{\bGamma}}$
in $\rep \hat{\bGamma}$.

We define the functions $S, T : \hat{\Gamma}_1 \to \{ \pm 1 \}$ by
\begin{align*}
S \beta & :=
\begin{cases}
S \alpha & \text{$\beta = \alpha [m]$ for $\alpha \in \Gamma_1$
and $m \in \bbZ$},
\\
- T \sigma & \text{$\beta = \sigma^* [m]$ for $\sigma \in \Sigma$
and $m \in \bbZ$},
\end{cases}
\\
\intertext{and} %
T \beta & :=
\begin{cases}
T \alpha & \text{$\beta = \alpha [m]$ for $\alpha \in \Gamma_1$
and $m \in \bbZ$},
\\
- S \sigma & \text{$\beta = \sigma^* [m]$ for $\sigma \in \Sigma$
and $m \in \bbZ$},
\end{cases}
\end{align*}
for $\beta \in \hat{\Gamma}_1$. One easily checks that
$(\hat{\Gamma}, \ul{\hat{R}})$, where $\ul{\hat{R}} := \bbZ R \cup
\{ \lambda \mid \lambda \in \Lambda \}$, is an almost gentle
quiver with string functions $S$ and $T$. By a path in
$\hat{\bGamma}$ we mean a path $(\hat{\Gamma}, \ul{\hat{R}})$.
Similarly, by a string (band) in $\hat{\bGamma}$ we mean a string
(band, respectively) in $(\hat{\Gamma}, \ul{\hat{R}})$. For a
string $\omega$ in $\bGamma$ and $m \in \bbZ$ we define the string
$\omega [m]$ in $\hat{\bGamma}$ in the obvious way.

With a string $\zeta$ in $\hat{\bGamma}$ we associate the
representation $V_\zeta$ of $\hat{\bGamma}$ in the following way.
First, we put $I_x := \{ i \in [0, \ell (\zeta)] \mid s
\zeta_{[i]} = x \}$ for $x \in \Gamma_0$. Then we define $V_\zeta$
by $V_\zeta (x) := k^{I_x}$ for $x \in \Gamma_0$ and
\[
(V_\zeta (\alpha))_{i, j} :=
\begin{cases}
\Id & \text{$i = j - 1$ and $\alpha = \alpha_j (\zeta)$, or}
\\
& \text{$i = j + 1$ and $\alpha = \alpha_{j + 1}^{-1} (\zeta)$},
\\
0 & \text{otherwise},
\end{cases}
\]
for $\alpha \in \Gamma_1$, $j \in I_{s \alpha}$, and $i \in I_{t
\alpha}$. If $i \in I_x$ for $x \in \Gamma_0$, then we denote by
$e_i (\zeta)$ the corresponding basis element of $V_\zeta (x)$.
Similarly, if $i \in I_x$ for $x \in \Gamma_0$, then $e_i^*
(\zeta)$ denotes the corresponding element of the basis of
$(V_\zeta (x))^*$ dual to $(e_j (\zeta))_{j \in I_x}$. The
representations of $\hat{\bGamma}$ of the above form are called
the string representations of $\hat{\bGamma}$.

For a string $\zeta$ in $\hat{\bGamma}$ of lenght $l$ we denote by
$\upsilon_\zeta$ the map $\upsilon : X_\zeta \to X_{\zeta^{-1}}$
given by $\upsilon (e_i (\zeta)) := e_{l - i} (\zeta^{-1})$ for $i
\in [0, l]$. Observe that $\upsilon_\zeta$ is an isomorphism for
each string $\zeta$ in $\hat{\bGamma}$ (and the inverse map is
given by $\upsilon_{\zeta^{-1}}$). Moreover, if $\zeta'$ and
$\zeta''$ are strings in $\hat{\bGamma}$ such that $V_{\zeta'}
\simeq V_{\zeta''}$, then either $\zeta' = \zeta''$ or
$\zeta'^{-1} = \zeta''$.

Let $\zeta'$ and $\zeta''$ be strings in $\hat{\bGamma}$ such that
$t \zeta' = t \zeta''$ and $T \zeta' = T \zeta''$. Put $l' := \ell
(\zeta')$, $l'' := \ell (\zeta'')$, and
\[
l := \max \{ i \in [0, \min (l', l'')] \mid \text{$\alpha_j
(\zeta') = \alpha_j (\zeta'')$ for each $j \in [1, i]$} \}.
\]
If $\zeta' \neq \zeta''$, then we write $\zeta' <_{\rt} \zeta''$
if either $l'' > l$ and $\alpha_{l + 1} (\zeta'') \in \Gamma_1$ or
$l' > l$ and $\alpha_{l + 1}^{-1} (\zeta') \in \Gamma_1$.
Moreover, we write $\zeta' \leq_{\rt} \zeta''$ if either $\zeta' =
\zeta''$ or $\zeta' <_{\rt} \zeta''$. If $\zeta' \leq_{\rt}
\zeta''$, then by $f_{\zeta', \zeta''}$ we denote the map $f :
V_{\zeta'} \to V_{\zeta''}$ given by $f (e_i (\zeta')) := e_i
(\zeta'')$ for $i \in [0, l]$ and $f (e_i (\zeta')) := 0$ for $i
\in [l + 1, l']$.

Dually, let $\zeta'$ and $\zeta''$ be strings in $\hat{\bGamma}$
such that $s \zeta' = s \zeta''$ and $S \zeta' = S \zeta''$. We
write $\zeta' \leq_{\rs} \zeta''$ if $\zeta'^{-1} \leq_{\rt}
\zeta''^{-1}$. If $\zeta' \leq_{\rs} \zeta''$, then we put
$g_{\zeta', \zeta''} := \upsilon_{\zeta''^{-1}} \circ
f_{\zeta'^{-1}, \zeta''^{-1}} \circ \upsilon_{\zeta'}$.

Let $\omega$ be a band in $\hat{\bGamma}$ and $\mu$ an
indecomposable automorphism of a finite dimensional vector space
$K$. We define the representation $W_{\omega, \mu}$ of
$\hat{\bGamma}$ as follows. First, for $x \in \Gamma_0$ we put
$J_x := \{ i \in [1, \ell (\omega)] \mid s \alpha_i (\omega) = x
\}$. Then $W_{\omega, \mu} (x) := K^{J_x}$ for $x \in Q_0$ and
\[
(W_{\omega, \mu} (\alpha))_{i, j} :=
\begin{cases}
\Id & \text{$j \in [2, \ell (\omega)]$, $i = j - 1$, and $\alpha =
\alpha_j (\omega)$, or}
\\
& \text{$j \in [1, \ell (\omega) - 1]$, $i = j + 1$, and $\alpha =
\alpha_{j + 1}^{-1} (\omega)$},
\\
\mu & \text{$j = 1$, $i = \ell (\omega)$, and $\alpha = \alpha_1
(\omega)$, or}
\\
& \text{$j = \ell (\omega)$, $i = 1$, and $\alpha = \alpha_1^{-1}
(\omega)$},
\\
0 & \text{otherwise},
\end{cases}
\]
for $\alpha \in \Gamma_1$, $j \in J_{s \alpha}$, and $i \in J_{t
\alpha}$. The representations of $\hat{\bGamma}$ of the above form
are called the band representations of $\hat{\bGamma}$.

Let $\Xi$ be the set consisting of all paths in $\hat{\bGamma}$
and a chosen full path starting at $y$ for each $y \in
\hat{\Gamma}_0$. For $\xi \in \Xi$ we denote by $\xi^*$ unique
$\xi' \in \Xi$ such that $t \xi' = s \xi$ and $\xi \xi'$ is a full
path. Observe that $(\sigma [m])^* = \sigma^* [m]$ for all $\sigma
\in \Sigma$ and $m \in \bbZ$. Moreover, $(\xi^*)^* = \nu \xi$ for
each $\xi \in \Xi$.

Fix $y \in \hat{\Gamma}_0$. If $y' \in \hat{\Gamma}_0$ and $\Xi
(y', y) := \{ \xi \in \Xi \mid s \xi = y', \; t \xi = y \}$, then
(the residue classes of) $(\xi)_{\xi \in \Xi (y', y)}$ form a
basis of $k \hat{\bGamma} (y', y)$. For each $y' \in
\hat{\Gamma}_0$ we identify $(\xi^*)_{\xi \in \Xi (y', y)}$ with
the basis of $(k \hat{\bGamma} (y', y))^*$ dual to $(\xi)_{\xi \in
\Xi (y', y)}$. This identification induces isomorphisms $(k
\hat{\bGamma} (y', y))^* \simeq k \hat{\bGamma} (\nu y, y')$, $y'
\in \hat{\Gamma}_0$, which extend to an isomorphism $Q_y \simeq
P_{\nu y}$, which we also treat as identification.

Let $\zeta$ be a string in $\hat{\bGamma}$ of positive length.
Then we have a presentation $\zeta = \xi_1 \xi_2 \cdots \xi_L$, $L
\in \bbN_+$, where $\xi_i$ is a directed string in $\hat{\bGamma}$
of positive length for each $i \in [1, L]$, and $\xi_i \xi_{i +
1}$ is not a directed string in $\hat{\bGamma}$ for each $i \in
[1, L - 1]$. In the above situation we put $L (\zeta) := L$, and
$\xi_i (\zeta) := \xi_i$ and $\xi_i^{-1} (\zeta) := \xi_i^{-1}$
for $i \in [1, L]$. Moreover, if $\zeta$ is a band in
$\hat{\bGamma}$, then we put $\zeta^{(0)} := \zeta$ and
$\zeta^{(i)} := \xi_{i + 1} \cdots \xi_L \xi_1 \cdots \xi_i$ for
$i \in [1, L - 1]$. We also put $L (\bone_{x, \varepsilon}) := 0$
for $x \in \hat{\Gamma}_0$ and $\varepsilon \in \{ \pm 1 \}$.

Now we define the operation $(-)^\times$ on the strings in
$\hat{\bGamma}$. First, we put $(\bone_{y, \varepsilon})^\times :=
\bone_{\nu y, - \varepsilon}$ for $y \in \hat{\Gamma}_0$ and
$\varepsilon \in \{ \pm 1 \}$. Next, if $\xi$ is a directed string
of positive length, then we put $\xi^\times := (\xi^*)^{-1}$ if
$\xi$ is a path in $\hat{\bGamma}$, and $\xi^\times :=
(\xi^{-1})^*$ if $\xi^{-1}$ is a path in $\hat{\bGamma}$. Finally,
if $\zeta$ is an arbitrary string in $\hat{\bGamma}$ of positive
length, then $\zeta^\times := \xi_1 (\zeta)^\times \cdots \xi_{L
(\zeta)} (\zeta)^\times$. The operation $(-)^\times$ is clearly
invertible. We denote the inverse operation by $(-)^+$.

We also need some additional operations on the strings in
$\hat{\bGamma}$. Let $\zeta$ be a string in $\hat{\bGamma}$, $l :=
\ell (\zeta)$, and $L := L (\zeta)$. If $L > 0$, then we denote by
$\partial \zeta$ the unique string $\zeta'$ in $\hat{\bGamma}$
such that $\zeta = \xi_1 (\zeta) \cdot \zeta'$. We put
\begin{align*}
\partial' \zeta & :=
\begin{cases}
\partial \zeta & \text{$L > 0$ and $\xi_1 (\zeta)$ is a path
in $\bGamma$},
\\
\zeta & \text{otherwise},
\end{cases}
\\
\intertext{and} %
\partial'' \zeta & :=
\begin{cases}
\partial \zeta & \text{$L > 0$ and $\xi_1^{-1} (\zeta)$ is a
path in $\bGamma$},
\\
\zeta & \text{otherwise}.
\end{cases}
\end{align*}
Next, we put
\begin{align*}
\delta_{\rt}' \zeta & :=
\begin{cases}
{}_{[1]} \zeta & \text{$l > 0$ and $\alpha_1 (\zeta) \in
\hat{\Gamma}_1$},
\\
\sigma_{t \zeta, - T \zeta} \zeta & \text{otherwise}.
\end{cases}
\\
\delta_{\rs}' \zeta & :=
\begin{cases}
\zeta_{[l - 1]} & \text{$l > 0$ and $\alpha_{\ell (\zeta)}^{-1}
(\zeta) \in \hat{\Gamma}_1$},
\\
\zeta \sigma_{s \zeta, - S \zeta}^{-1} & \text{otherwise},
\end{cases}
\\
\delta_{\rt}'' \zeta & :=
\begin{cases}
{}_{[1]} \zeta & \text{$l > 0$ and $\alpha_1^{-1} (\zeta) \in
\hat{\Gamma}_1$},
\\
\sigma_{t \zeta, - T \zeta}'^{-1} \zeta & \text{otherwise}.
\end{cases}
\\
\intertext{and} %
\delta_{\rs}'' \zeta & :=
\begin{cases}
\zeta_{[l - 1]} & \text{$l > 0$ and $\alpha_{\ell (\zeta)} (\zeta)
\in \hat{\Gamma}_1$},
\\
\zeta \sigma_{s \zeta, - S \zeta}' & \text{otherwise}.
\end{cases}
\end{align*}
Moreover, we put $\delta' \zeta := \delta_{\rs}' \delta_{\rt}'
\zeta = \delta_{\rt}' \delta_{\rs}' \zeta$ and $\delta'' \zeta :=
\delta_{\rs}'' \delta_{\rt}'' \zeta = \delta_{\rt}''
\delta_{\rs}'' \zeta$. Finally, we put $\Delta \zeta := \delta'
(\zeta^\times)$. Observe that $\Delta$ is invertible and
$\Delta^{-1} \zeta := \delta'' (\zeta^+)$.

Let $\zeta$ be a string in $\hat{\bGamma}$. We describe a
projective cover $\pi_\zeta : P_\zeta \to V_\zeta$ of $V_\zeta$.
Write $\zeta = \xi_1 \xi_2^{-1} \cdots \xi_{2 L - 1} \xi_{2
L}^{-1}$, where $L \in \bbN_+$ and $\xi_1$, \ldots, $\xi_{2 L}$
are simple strings in $\bGamma$ such that $\ell (\xi_i) > 0$ for
each $i \in [2, 2 L - 1]$. Let $l_i := \sum_{j \in [1, 2 i - 1]}
\ell (\xi_j)$ for $i \in [1, L]$. Then $\pi_\zeta : P_\zeta \to
V_\zeta$, where $P_\zeta := \bigoplus_{i \in [1, L]} P_{s \xi_{2 i
- 1}}$ and $(\pi_\zeta)_i$ corresponds to $e_{l_i} (\zeta)$ under
the canonical isomorphism $\Hom_{\bGamma} (P_{s \xi_{2 i - 1}},
V_\zeta) \simeq V_\zeta (s \xi_{2 i - 1})$ for each $i \in [1,
L]$, is the minimal projective cover of $V_\zeta$. Observe that
$\Omega V_\zeta := \Ker \pi_\zeta \simeq V_{\Delta^{-1} \zeta}$.
We identify $\Omega V_\zeta$ with $V_{\Delta^{-1} \zeta}$.
Further, $P_\zeta$ is an injective envelope of $V_{\Delta^{-1}
\zeta}$. More precisely, if $\Delta^{-1} \zeta = \xi_1'^{-1}
\xi_2' \cdots \xi_{2 L - 1}'^{-1} \xi_{2 L}'$ for $L \in \bbN_+$
and simple strings $\xi_1'$, \ldots, $\xi_{2 L}'$ in
$\hat{\bGamma}$ such that $\ell (\xi_i') > 0$ for each $i \in [2,
2 L - 1]$, then $P_\zeta = \bigoplus_{i \in [1, L]} Q_{t \xi_{2 i
- 1}'}$. Moreover, if $l_i' := \sum_{j \in [1, 2 i - 1]} \ell
(\xi_j')$ for $i \in [1, L]$ and $\iota_\zeta : V_{\Delta^{-1}
\zeta} \to P_\zeta$ is such that $(\iota_\zeta)_i$ corresponds to
$(-1)^i e_{l_i'}^* (\Delta^{-1} \zeta)$ under the canonical
isomorphism $(V_{\Delta^{-1} \zeta} (t \xi_{2 i - 1}))^* \simeq
\Hom_{\hat{\bGamma}} (V_{\Delta^{-1} \zeta}, Q_{t \xi_{2 i - 1}})$
for each $i \in [1, L]$, then the sequence $0 \to V_{\Delta^{-1}
\zeta} \xrightarrow{\iota_\zeta} P_\zeta \xrightarrow{\zeta}
V_\zeta \to 0$ is exact. We will use sequences of the above form
to calculate the action of $\Omega$ on morphisms in $\urep
\hat{\bGamma}$. In particular, it follows that $\Omega f_{\xi',
\xi''}$ and $f_{\Delta^{-1} \xi', \Delta^{-1} \xi''}$ ($\Omega
g_{\xi', \xi''}$ and $g_{\Delta^{-1} \xi', \Delta^{-1} \xi''}$)
coincide up to sign for strings $\xi'$ and $\xi''$ in
$\hat{\bGamma}$ such that $t \xi' = t \xi''$, $T \xi' = T \xi''$,
and $\xi' \leq_{\rt} \xi''$ ($s \xi' = s \xi''$, $S \xi' = S
\xi''$, and $\xi' \leq_{\rs} \xi''$).

Similarly as above we show that $\Omega W_{\zeta, \mu} \simeq
W_{\zeta^+, (-1)^{L (\zeta) / 2} \mu^{-1}}$.

\section{The Happel functor} \label{sect_Happel}

Throughout this section we fix a gentle bound quiver $\bGamma  =
(\Gamma, R)$ together with string functions $S$ and $T$. We also
denote by $\Sigma$ the set of maximal paths in $\bGamma$.

Let $\calD^b (\bGamma)$ denote the derived category of $\rep
\bGamma$. It is known that $\calD^b (\bGamma)$ is a triangulated
category and $\calK^b (\bGamma)$ can be viewed as a full
triangulated subcategory of $\calD^b (\bGamma)$. We identify $M
\in \rep \bGamma$ with the complex in $\calD^b (\bGamma)$
concentrated in degree $0$. In~\cite{Happel1988} Happel
constructed a fully faithful triangle functor $\calD^b (\bGamma)
\to \urep \hat{\bGamma}$ which extends the inclusion functor $\rep
\bGamma \to \urep \hat{\bGamma}$. By $\Psi$ we denote the
restriction of this functor to $\calK^b (\bGamma)$.

Let $\omega$ be a homotopy string $\bGamma$. We define the string
$\psi \omega$ in $\hat{\bGamma}$ by induction on $L (\omega)$ as
follows. First, $\psi \omega := \sigma_\omega [0] \cdot
(\sigma_{\omega^{-1}})^{-1} [0]$ if $L (\omega) = 0$. Next, if $L
(\omega) > 0$, then
\[
\psi \omega :=
\begin{cases}
\sigma_\omega [0] \cdot (\sigma [0])^+ \cdot \delta_{\rs}''
(\zeta^+) & \text{$\sigma$ is a path in $\bGamma$},
\\
\sigma_\omega [0] \cdot \sigma [0] \cdot \delta_{\rs}'
(\zeta^\times) & \text{$\sigma^{-1}$ is a path in $\bGamma$ and
$\ell (\zeta) > 0$},
\\
\sigma_\omega [0] \cdot \sigma_{[\ell (\sigma) - 1]} [0] &
\text{$\sigma^{-1}$ is a path in $\bGamma$ and $\ell (\zeta) =
0$},
\end{cases}
\]
where $\sigma := \sigma_1 (\omega)$ and $\zeta := \partial' (\psi
({}^{[1]} \omega))$.

The meaning of the above assignment is explained in the following.

\begin{proposition} \label{proposition_psi}
Let $\omega$ be a homotopy string in $\bGamma$ and $m \in \bbZ$.
Then $\Psi X_{m, \omega} \simeq V_{\Delta^{-m} (\psi \omega)}$.
Moreover, under the isomorphisms of the above form, we have the
following:
\begin{enumerate}

\item \label{prop_psi_one}
If $\sigma$ is a path in $\bGamma$ of positive length such that
the composition $\sigma \omega$ is defined, then $\Psi F'_{m,
\sigma, \omega} = \varepsilon f_{\Delta^{-m} (\psi \bone_{t
\sigma, - T \sigma}), \Delta^{-m} (\psi \omega)}$ for some
$\varepsilon \in \{ \pm 1 \}$.

\item
If $\sigma$ is a path in $\bGamma$ of positive length such that
$\sigma^{-1} \omega$ is defined, then $\Psi F''_{m, \sigma,
\omega} = \varepsilon f_{\Delta^{-m} (\psi \omega), \Delta^{-m}
(\psi \bone_{s \sigma, S \sigma})}$ for some $\varepsilon \in \{
\pm 1 \}$.

\item \label{prop_psi_four}
If $\sigma$ is a path in $\bGamma$ of positive length such that
$\omega \sigma^{-1}$ is defined, then $\Psi G_{m, \sigma, \omega}'
= \varepsilon g_{\Delta^{- \deg \omega - m} (\psi \bone_{t \sigma,
T \sigma}), \Delta^{-m} (\psi \omega)}$ for some $\varepsilon \in
\{ \pm 1 \}$.

\item
If $\sigma$ is a path in $\bGamma$ of positive length such that
$\omega \sigma$ is defined, then $\Psi G_{m, \sigma, \omega}'' =
\varepsilon g_{\Delta^{-m} (\psi \omega), \Delta^{- \deg \omega -
m} (\psi \bone_{s \sigma, - S \sigma})}$ for some $\varepsilon \in
\{ \pm 1 \}$.

\end{enumerate}
\end{proposition}

\begin{proof}
We prove the above claims by induction on $L (\omega)$. If $L
(\omega) = 0$, then the proofs are immediate (observe that $X_{m,
\omega} \simeq X_{0, \omega} [-m]$ for each $m \in \bbZ$). It also
follows easily that $\Psi \Upsilon_{0, \omega}$ and
$\upsilon_{\psi \omega}$ coincide in this case.

Now assume that $L (\omega) > 0$ and let $\sigma' := \sigma_1
(\omega)$ and $\omega' := {}^{[1]} \omega$. We also assume that
$\sigma'$ is a path in $\bGamma$, the case when $\sigma'^{-1}$ is
a path in $\bGamma$ is similar.

By induction hypothesis $\Psi X_{0, \bone_{t \sigma', - T
\sigma'}} \simeq V_{\psi \bone_{t \sigma', - T \sigma'}}$ and
$\Psi X_{0, \omega'} \simeq V_{\psi \omega'}$. Moreover, if $f :=
\Psi F_{0, \sigma', \omega'}'$, then $f$ equals up to sign
$f_{\psi \bone_{t \sigma', - T \sigma'}, \psi \omega'}$. Let
$\zeta := \Delta (\psi \bone_{t \sigma', - T \sigma'})$ and
\[
\SelectTips{cm}{} %
\vcenter{\xymatrix{V_{\psi \bone_{t \sigma', - T \sigma'}}
\ar[r]^{\iota_\zeta} \ar[d]^f & P_\zeta \ar[r]^{\pi_\zeta}
\ar[d]^g & V_\zeta \ar[d]^{\Id} \\ V_{\psi \omega'} \ar[r]^{f'} &
V \ar[r]^{f''} & V_\zeta}}
\]
be the push-out diagram. Then
\[
V_{\psi \bone_{t \sigma', - T \sigma'}} \xrightarrow{f} V_{\psi
\omega'} \xrightarrow{f'} V \xrightarrow{f''} V_{\zeta}
\]
is a triangle in $\urep \hat{\bGamma}$. Direct calculations show
that $V \simeq V_{\Delta (\psi \omega)}$ in $\urep \hat{\bGamma}$.
Moreover, by choosing $g$ in an appropriate way we get that $f'$
and $f''$ equal up to sign $g_{\psi \omega', \Delta (\psi
\omega)}$ and $\upsilon_{\zeta^{-1}} \circ f_{\Delta (\psi
\omega), \zeta^{-1}}$, respectively. Since we have an isomorphism
$X_{m, \omega} \simeq X_{-1, \omega} [- m - 1]$ for each $m \in
\bbZ$ and a triangle
\[
X_{0, \bone_{t \sigma', - T \sigma'}} \xrightarrow{F_{0, \sigma',
\omega'}'} X_{0, \omega'} \xrightarrow{F'} X_{-1, \omega}
\xrightarrow{F''} X_{-1, \bone_{t \sigma', - T \sigma'}}
\]
in $\calK^b (\bGamma)$, we get the first claim. Moreover, under
the appropriate isomorphisms, $\Psi F' = f'$ and $\Psi F'' = f''$.

Let $\sigma$ be a path in $\bGamma$ of positive length such that
$s \sigma = t \omega$ and $S \sigma = T \omega$. Then we have the
commutative diagram
\[
\xymatrix{X_{0, \bone_{t \sigma, - T \sigma}} \ar[r] \ar[d]^{F_{0,
\sigma, \bone_{t \sigma', - T \sigma'}}'} & 0 \ar[r] \ar[d] &
X_{-1, \bone_{t \sigma, - T \sigma}} \ar[r]^{\Id} \ar[d]^{F_{0,
\sigma, \omega}'} & X_{-1, \bone_{t \sigma, - T \sigma}}
\ar[d]^{F_{-1, \sigma, \bone_{t \sigma', - T \sigma'}}'} \\ X_{0,
\bone_{t \sigma', - T \sigma'}} \ar[r]^{F_{0, \sigma', \omega'}'}
& X_{0, \omega'} \ar[r]^{F'} & X_{-1, \omega} \ar[r]^{F''} &
X_{-1, \bone_{t \sigma', - T \sigma'}}}
\]
in $\calK^b (\bGamma)$, which gives rise to the commutative
diagram
\[
\xymatrix{V_{\psi \bone_{t \sigma, - T \sigma}} \ar[r]
\ar[d]^{\Psi F_{0, \sigma, \bone_{t \sigma', - T \sigma'}}'} & 0
\ar[r] \ar[d] & V_{\Delta (\psi \bone_{t \sigma, - T \sigma})}
\ar[r]^{\Id} \ar[d]^{\Psi F_{0, \sigma, \omega}'} & V_{\Delta
(\psi \bone_{t \sigma, - T \sigma})} \ar[d]^{\Psi F_{-1, \sigma,
\bone_{t \sigma', - T \sigma'}}'} \\ V_{\psi \bone_{t \sigma', - T
\sigma'}} \ar[r]^f & V_{\psi \omega'} \ar[r]^{f'} & V_{\Delta
(\psi \omega)} \ar[r]^{f''} & V_\zeta}
\]
in $\urep \hat{\bGamma}$. By induction hypothesis $\Psi F_{0,
\sigma, \bone_{t \sigma', - T \sigma'}}'$ and $\Psi F_{-1, \sigma,
\bone_{t \sigma', - T \sigma'}}'$ equal up to sign $f_{\psi
\bone_{t \sigma, - T \sigma}, \psi \bone_{t \sigma', - T
\sigma'}}$ and $f_{\Delta (\psi \bone_{t \sigma, - T \sigma}),
\zeta}$, respectively. Moreover, $\Psi F_{m, \sigma, \bone_{t
\sigma', - T \sigma'}}' = \Omega^{m + 1} \Psi F_{-1, \sigma,
\bone_{t \sigma', - T \sigma'}}'$ and the above diagram determines
$\Psi F_{-1, \sigma, \bone_{t \sigma', - T \sigma'}}'$,
hence~\eqref{prop_psi_one} follows. The remaining claims are
proved similarly.
\end{proof}

Observe that $\omega' = \omega''$ if $\psi \omega' = \psi
\omega''$ for homotopy strings $\omega'$ and $\omega''$ in
$\bGamma$. In fact we have more.

\begin{lemma}
Let $\omega'$ and $\omega''$ be homotopy strings in $\bGamma$ and
$m \in \bbZ$. If $\psi \omega' = \Delta^m (\psi \omega'')$, then
$m = 0$ and $\omega' = \omega''$.
\end{lemma}

\begin{proof}
In light of the preceding remark it is sufficient to prove that $m
= 0$. Without loss of generality we may assume that $m \leq 0$.
Observe that $t (\psi \omega') \in \Gamma_0 [0]$. Moreover, easy
induction shows that $t (\Delta^m (\psi \omega'')) \in \Gamma_0
[-m'']$ for some $m'' \in \bbN_+ 0$ provided $m < 0$, hence the
claim follows.
\end{proof}

\begin{corollary} \label{coro_phiomegainverse}
Let $\omega$ be a homotopy string in $\bGamma$. Then $(\psi
\omega)^{-1} = \Delta^{- \deg \omega} (\psi \omega^{-1})$.
\end{corollary}

\begin{proof}
Recall that $X_{0, \omega} \simeq X_{\deg \omega, \omega^{-1}}$,
thus Proposition~\ref{proposition_psi} implies that $V_{\psi
\omega} \simeq V_{\Delta^{- \deg \omega} (\psi \omega^{-1})}$.
Hence, either $\psi \omega = \Delta^{- \deg \omega} (\psi
\omega^{-1})$ or $(\psi \omega)^{-1} = \Delta^{- \deg \omega}
(\psi \omega^{-1})$. In the former case the previous lemma implies
that $\omega^{-1} = \omega$, which is impossible.
\end{proof}

Now we calculate the images of the band complexes. For this we
need an additional function $\psi'$ between the homotopy strings
in $\bGamma$ of positive length and the strings in
$\hat{\bGamma}$. Let $\omega$ be a homotopy string in $\bGamma$ of
positive length. Put $L := L (\omega)$ and $\sigma := \sigma_1
(\omega)$. If $L = 1$, then $\psi' \omega := (\sigma [0])^+$
provided $\sigma$ is a path in $\bGamma$, and $\psi' \omega :=
\sigma [0]$ provided $\sigma^{-1}$ is a path in $\bGamma$. If $L >
1$, then
\[
\psi' \omega :=
\begin{cases}
\psi' \sigma \cdot (\psi' ({}^{[1]} \omega))^+ & \text{$\sigma$ is
a path in $\bGamma$},
\\
\psi' \sigma \cdot (\psi' ({}^{[1]} \omega))^\times &
\text{$\sigma^{-1}$ is a path in $\bGamma$}.
\end{cases}
\]
One can easily deduce the formula for $\psi' \omega$ as follows:
\[
\psi' \omega = (\sigma_1 (\omega))^{\times n_1} \cdots (\sigma_L
(\omega))^{\times n_L},
\]
where
\[
n_i :=
\begin{cases}
- \deg \omega^{[i]} & \text{$\sigma_i (\omega)$ is a path in
$\bGamma$},
\\
- \deg \omega^{[i]} - 1 & \text{$\sigma_i^{-1} (\omega)$ is a path
in $\bGamma$},
\end{cases}
\]
for $i \in [1, L]$, and $(-)^{\times n}$ denotes the $n$-th power
of $(-)^\times$ for $n \in \bbZ$. The above formula implies
immediately that $\psi' \omega^{-1} = (\psi' \omega)^{-1}$ if
$\deg \omega = 0$. Moreover, if $\deg \omega = 0$ and $\omega'$ is
a homotopy string in $\bGamma$ of positive length such that the
composition $\omega \omega'$ is defined, then $\psi' (\omega
\omega') = \psi' \omega \cdot \psi' \omega'$. Finally, it follows
that $\psi' \omega$ is a band in $\hat{\bGamma}$ if $\omega$ is a
homotopy band in $\bGamma$.

We also need the following property of $\psi'$.

\begin{lemma}
Let $\omega$ be a homotopy string in $\bGamma$ such that $L
(\omega) > 0$ and $\deg \omega = 0$. If $\omega'$ is a homotopy
string in $\bGamma$ such that the composition $\omega \omega'$ is
defined, then $\psi (\omega \omega') = \sigma_\omega [0] \cdot
\psi' \omega \cdot \partial' (\psi \omega')$. In particular, $\psi
\omega = \sigma_\omega [0] \cdot \psi' \omega \cdot
(\sigma_{\omega^{-1}})^{-1} [0]$.
\end{lemma}

\begin{proof}
If $L (\omega) = 2$, then the claim follows by direct
calculations. Observe that in this case either $\omega = \sigma'
\sigma''^{-1}$ or $\omega = \sigma'^{-1} \sigma''$ for paths
$\sigma'$ and $\sigma''$ in $\bGamma$ of positive length.

Now assume that $L := L (\omega) > 2$. There are some cases to
consider.

First assume that there exists $i \in [2, L - 2]$ such that $\deg
\omega^{[i]} = 0$. Observe that $\deg {}^{[i]} \omega = 0$. Then
we have the following sequence of equalities
\begin{align*}
\psi (\omega \omega') & \induction \sigma_\omega [0] \cdot \psi'
\omega^{[i]} \cdot \partial' (\psi' ({}^{[i]} \omega \omega'))
\induction
\\
& \induction \sigma_\omega [0] \cdot \psi' \omega^{[i]} \cdot
\psi' {}^{[i]} \omega \cdot  \partial' (\psi \omega') =
\sigma_\omega [0] \cdot \psi' \omega \cdot \partial' (\psi'
\omega')
\end{align*}
hence the claim follows in this case.

If the above condition is not satisfied, then $\deg {}^{[1]}
\omega^{[L - 1]} = 0$. Put $\sigma' := \sigma_1 (\omega)$ and
$\sigma'' := \sigma_{L} (\omega)$. Assume in addition that
$\sigma'$ is a path in $\bGamma$, the other case is similar. Then
$\sigma''^{-1}$ is a path in $\bGamma$. Moreover, if $\ell
(\partial' (\psi \omega')) > 0$, then we have the following
sequence of equalities
\begin{align*}
\psi (\omega \omega') & \definition \sigma_{\omega} [0] \cdot
(\sigma' [0])^+ \cdot \delta_{\rs}'' ((\partial' (\psi ({}^{[1]}
\omega \cdot \omega')))^+)
\\
& \induction \sigma_{\omega} [0] \cdot (\sigma' [0])^+ \cdot
\delta_{\rs}'' ((\psi' ({}^{[1]} \omega^{[L - 1]}))^+ \cdot
(\partial' (\psi (\sigma'' \omega')))^+)
\\
& \definition \sigma_{\omega} [0] \cdot \psi' \sigma' \cdot
\delta_{\rs}'' ((\psi' ({}^{[1]} \omega^{[L - 1]}))^+ \cdot
(\sigma'' [0])^+ \cdot (\delta_{\rs}' ((\partial' (\psi
\omega'))^\times))^+)
\\
& = \sigma_{\omega} [0] \cdot \psi' \sigma' \cdot (\psi' ({}^{[1]}
\omega^{[L - 1]}))^+ \cdot (\psi' \sigma'')^+ \cdot
\partial'(\psi \omega')
\\
&  = \sigma_{\omega} [0] \cdot \psi' \sigma' \cdot (\psi'
({}^{[1]} \omega))^+ \cdot \partial' (\psi \omega') \definition
\sigma_{\omega} [0] \cdot \psi' \omega \cdot \partial' (\psi'
\omega').
\end{align*}
On the other hand, if $\ell (\partial' (\psi' \omega')) = 0$, then
by repeating some of the calculations above we get
\begin{align*}
\psi (\omega \omega') & = \sigma_{\omega} [0] \cdot \psi' \sigma'
\cdot \delta_{\rs}'' ((\psi' ({}^{[1]} \omega^{[L - 1]}))^+ \cdot
(\sigma_{[\ell (\sigma'') - 1]}'' [0])^+)
\\
& = \sigma_{\omega} [0] \cdot \psi' \sigma' \cdot (\psi' ({}^{[1]}
\omega^{[L - 1]}))^+ \cdot (\sigma'' [0])^+
\\
& = \sigma_{\omega} [0] \cdot \psi' \omega \cdot \partial' (\psi'
\omega'),
\end{align*}
what finishes the proof.
\end{proof}

The above lemma will be used in the proof of the following.

\begin{proposition} \label{prop_band}
Let $\omega$ be a homotopy band in $\bGamma$ and $\mu$ an
indecomposable automorphism of a finite dimensional vector space.
Let
\[
\varepsilon :=
\begin{cases}
1 & \text{$\sigma_1 (\omega)$ is a path in $\bGamma$},
\\
-1 & \text{$\sigma_1^{-1} (\omega)$ is a path in $\bGamma$}.
\end{cases}
\]
Then $\Psi Y_{0, \omega, \mu} \simeq W_{\psi' \omega, \varepsilon'
\mu^{\varepsilon}}$ for some $\varepsilon' \in \{ \pm 1 \}$
depending only on $\omega$.
\end{proposition}

\begin{proof}
Put $L := L (\omega)$. Let $\sigma' := \sigma_1 (\omega)$,
$\sigma'' := \sigma_{L (\omega)}^{-1} (\omega)$, and $\omega' :=
{}^{[1]} \omega^{[L - 1]}$. Observe that $\deg \omega' = 0$.
Consequently, $\psi' \omega' = \sigma_{\omega'} [0] \cdot \psi'
\omega' \cdot (\sigma_{\omega'^{-1}})^{-1} [0]$ if $L > 2$ (by the
previous lemma) and $\psi' \omega' = \sigma_{\omega'} [0] \cdot
(\sigma_{\omega'^{-1}})^{-1} [0]$ otherwise (by definition). Let
$K$ be the domain of $\mu$.

We assume that $\sigma'$ is a path in $\bGamma$ --- the other case
is similar. Then $\sigma''$ is also a path in $\bGamma$ and we
have a triangle
\[
X_{0, \bone_{t \sigma', - T \sigma'}} \otimes_k K \xrightarrow{F}
X_{0, \omega'} \otimes_k K \to Y_{-1, \omega, \mu} \to X_{-1,
\bone_{t \sigma', - T \sigma'}} \otimes K
\]
in $\calK^b (\bGamma)$, where $F := F_{0, \sigma', \omega'}'
\otimes \mu + G_{0, \sigma'', \omega'} \otimes \Id$. Observe that
Proposition~\ref{proposition_psi} implies that $\Psi X_{0,
\bone_{t \sigma', - T \sigma'}} \simeq V_{\psi \bone_{t \sigma', -
T \sigma'}}$ and $\Psi X_{0, \omega'} \simeq V_{\psi \omega'}$. If
$f := \Psi F$, then under the above isomorphisms $f =
\varepsilon_1 f_{\psi \bone_{t \sigma', - T \sigma'}, \psi
\omega'} \otimes \mu + \varepsilon_2 g_{\psi \bone_{t \sigma', - T
\sigma'}, \psi \omega'} \otimes \Id$ for some $\varepsilon_1,
\varepsilon_2 \in \{ \pm 1 \}$ (depending only on $\omega$)
according to
Proposition~\ref{proposition_psi}~\eqref{prop_psi_one}
and~\eqref{prop_psi_four}. Since we have the triangle
\begin{multline*}
V_{\psi \bone_{t \sigma', - T \sigma'}} \otimes K \xrightarrow{f}
V_{\psi \omega'} \otimes K \to W_{\sigma' [0] \cdot \psi' \omega'
\cdot \sigma''^{-1} [0], \varepsilon_0 \mu^{-1}}
\\
\to \Omega^{-1} (V_{\psi \bone_{t \sigma', - T \sigma'}} \otimes
K)
\end{multline*}
in $\urep \hat{\bGamma}$, where $\varepsilon_0 := - \varepsilon_1
\varepsilon_2$, $\Psi Y_{-1, \omega, \mu} \simeq W_{\sigma' [0]
\cdot \psi' \omega' \cdot \omega''^{-1} [0], \varepsilon_0
\mu^{-1}}$ (in the calculation of this triangle we use the form of
$\psi' \omega'$ calculated at the beginning of the proof). Since
$Y_{0, \omega, \mu} \simeq \Psi Y_{-1, \omega, \mu} [-1]$ and
\[
\Omega W_{\sigma' [0] \cdot \psi' \omega' \cdot \sigma''^{-1} [0],
\varepsilon_0 \mu^{-1}} \simeq W_{\psi' \omega, (-1)^{L (\omega) /
2} \varepsilon_0 \mu},
\]
the claim follows.
\end{proof}

\begin{corollary} \label{coro_band}
Let $m \in \bbZ$, $\omega$ be a homotopy band in $\bGamma$, and
$\mu$ an automorphism of a finite dimensional vector space. Then
$\Psi Y_{m, \omega, \mu} \simeq W_{(\psi' \omega)^{\times (-m)},
\varepsilon' \mu^{\varepsilon''}}$ for some $\varepsilon',
\varepsilon'' \in \{ \pm 1 \}$ depending only on $\omega$ and $m$.
\end{corollary}

\begin{proof}
Since $Y_{m, \omega, \mu} \simeq Y_{0, \omega, \mu} [-m]$, the
claim follows from the previous proposition.
\end{proof}

\section{Almost split triangles} \label{sect_almostsplit}

Throughout this section we fix a gentle bound quiver $\bGamma =
(\Gamma, R)$ together with string functions $S$ and $T$. Moreover,
$\psi$ is the map which associates with a homotopy string in
$\bGamma$ a string in $\hat{\bGamma}$ as defined in
Proposition~\ref{proposition_psi}.

Our aim in this section is to determine the shape of the almost
split triangles in $\calK^b (\bGamma)$. In order to restrict the
number of technical definitions we will not describe the maps
appearing in the almost split triangles, however the interested
reader can easily check that ``natural'' candidates are the
correct ones.

We first recall that a triangle $X \to Y \to Z \to X [1]$ is
$\calK^b (\bGamma)$ is an almost split triangle in $\calK^b
(\bGamma)$ if and only if the corresponding triangle $\Psi X \to
\Psi Y \to \Psi Z \to \Omega^{-1} \Psi X$ in $\urep \hat{\bGamma}$
is an almost split triangle in $\urep
\hat{\bGamma}$~\cite{HappelKellerReiten2008}*{Proposition~5.2},
where similarly as in the previous section $\Psi$ denotes the
restriction of the Happel functor to $\calK^b (\bGamma)$.

As a consequence we immediately obtain the following.

\begin{maintheorem}[Part I: Band complexes]
Let $m \in \bbZ$ and $\omega$ be a homotopy band in $\bGamma$. If
$0 \to \mu' \to \bigoplus_{i \in [1, n]} \mu_i \to \mu'' \to 0$ is
an almost split sequence in the category of the automorphisms of
finite dimensional vector spaces, where $\mu'$, $\mu_1$, \ldots,
$\mu_n$ and $\mu''$ are indecomposable automorphisms of finite
dimensional vector spaces, then
\[
Y_{m, \omega, \mu'} \to \bigoplus_{i \in [1, n]} Y_{m, \omega,
\mu_i} \to Y_{m, \omega, \mu''} \to Y_{m, \omega, \mu'} [1]
\]
is an almost split triangle in $\calK^b (\bGamma)$.
\end{maintheorem}

\begin{proof}
It follows from the above remark and Corollary~\ref{coro_band},
since for each band $\zeta$ in $\hat{\bGamma}$ and $\varepsilon,
\varepsilon' \in \{ \pm 1 \}$ we have an almost split triangle
\[
W_{\zeta, \varepsilon' \mu'^\varepsilon} \to \bigoplus_{i \in [1,
n]} W_{\zeta, \varepsilon' \mu_i^\varepsilon} \to W_{\zeta,
\varepsilon' \mu''^\varepsilon} \to \Omega^{-1} W_{\zeta,
\varepsilon' \mu'^\varepsilon}
\]
in $\urep \hat{\bGamma}$
(compare~\citelist{\cite{ButlerRingel1987}
\cite{WaldWaschbusch1985}*{Theorem~4.1 (2)}}).
\end{proof}

It remains to describe the almost split triangles in $\calK^b
(\bGamma)$ involving the string complexes. We first describe the
almost split triangles $\urep \hat{\bGamma}$ involving the string
representations.

For a string $\zeta$ in $\hat{\bGamma}$ we put
\[
{}_+ \zeta :=
\begin{cases}
\sigma_{s \alpha', - S \alpha'} \alpha'^{-1} \zeta &
\text{$\alpha' \neq \varnothing$},
\\
{}_{[1]} (\partial'' \zeta) & \text{$\alpha' = \varnothing$ and
$\ell (\partial'' \zeta) > 0$},
\\
\varnothing & \text{$\alpha' = \varnothing$ and $\ell (\partial''
\zeta) = 0$},
\end{cases}
\]
where
\[
\alpha' :=
\begin{cases}
\alpha_{t \zeta, - T \zeta}' & \text{$\alpha_{t \zeta, - T \zeta}'
\neq \varnothing$ and $\alpha_{t \zeta, - T \zeta}'^{-1} \zeta$ is
a string in $\hat{\bGamma}$},
\\
\varnothing & \text{otherwise}.
\end{cases}
\]
Next, we put $\zeta_+ := ({}_+ (\zeta^{-1}))^{-1}$ for a string
$\zeta$ in $\hat{\bGamma}$. Finally, if $\zeta$ is a string in
$\hat{\bGamma}$, then
\[
{}_+ \zeta_+ :=
\begin{cases}
({}_+ \zeta)_+ & \text{${}_+ \zeta \neq \varnothing$},
\\
{}_+ (\zeta_+) & \text{$\zeta_+ \neq \varnothing$}.
\end{cases}
\]
We leave it to the reader to verify that the above definition is
correct and ${}_+ \zeta_+ \neq \varnothing$ for each string
$\zeta$ in $\hat{\bGamma}$. Moreover,
\citelist{\cite{ButlerRingel1987}
\cite{WaldWaschbusch1985}*{Theorem~4.1 (1)}} imply that for each
string $\zeta$ in $\hat{\bGamma}$ we have an almost split triangle
in $\urep \hat{\bGamma}$ of the form
\[
V_\zeta \to V_{{}_+ \zeta} \oplus V_{\zeta_+} \to V_{{}_+ \zeta_+}
\to \Omega^{-1} V_\zeta,
\]
where $V_\varnothing := 0$.

We translate the above construction to $\calK^b (\bGamma)$.

Let $\omega$ be a homotopy band in $\bGamma$. We denote by $r
(\omega)$ the maximal $i \in [0, \ell (\omega)]$ such that
$\alpha_j (\omega) \in \Gamma_1$ for each $j \in [1, i]$ and
$\alpha_j (\omega) \alpha_{j + 1} (\omega) \in R$ for each $j \in
[1, i - 1]$. We put
\[
\omega' :=
\begin{cases}
{}_{[1]} (\sigma_{r (\omega)} (\omega)) \cdot {}^{[r (\omega)]}
\omega & \text{$r (\omega) > 0$},
\\
\omega & \text{$r (\omega) = 0$},
\end{cases}
\]
and $\sigma := \sigma_\omega$. Then we define
\begin{gather*}
{}_+ \omega :=
\begin{cases}
\theta_{t \sigma, - T \sigma}^{-1} \sigma \omega & \text{$\ell
(\sigma) > 0$},
\\
\theta_{t \omega', - T \omega'}^{-1} \omega' & \text{$\ell
(\sigma) = 0$, $\ell (\theta_{t \omega', - T \omega'}) > 0$, and
$\ell (\omega') > 0$},
\\
({}^{[1]} \theta_{t \omega', - T \omega'})^{-1} & \text{$\ell
(\sigma) = 0$, $\ell (\theta_{t \omega', - T \omega'}) > 0$, and
$\ell (\omega') = 0$},
\\
{}^{[1]} \omega' & \text{$\ell (\sigma) = 0$, $\ell (\theta_{t
\omega', - T \omega'}) = 0$, $\ell (\omega') > 0$},
\\
& \quad \text{and $\alpha_1^{-1} (\omega') \in \Gamma_1$},
\\
\omega' & \text{$\ell (\sigma) = 0$, $\ell (\theta_{t \omega', - T
\omega'}) = 0$, $\ell (\omega') > 0$},
\\
& \quad \text{and $\alpha_1 (\omega') \in \Gamma_1$},
\\
\varnothing & \text{$\ell (\sigma) = 0$, $\ell (\theta_{t \omega',
- T \omega'}) = 0$, and $\ell (\omega') = 0$}.
\end{cases}
\\
\intertext{and} %
m' (\omega) :=
\begin{cases}
\ell (\theta_{t \sigma, - T \sigma}) - 1 & \text{$\ell (\sigma) >
0$},
\\
\ell (\theta_{t \omega', - T \omega'}) + r (\omega) - 1 &
\text{$\ell (\sigma) = 0$}.
\end{cases}
\end{gather*}

We first prove that the above definitions are correct.

\begin{lemma}
Let $x \in \Gamma_0$ and $\varepsilon \in \{ \pm 1 \}$. If
$\alpha_{x, \varepsilon} = \varnothing$, then $\theta_{x,
\varepsilon} \neq \varnothing$.
\end{lemma}

\begin{proof}
We show that $\ell (\theta) \leq |\Gamma_1|$ for each $\theta \in
\Theta_{x, \varepsilon}$. Assume this is not the case and fix
$\theta \in \Theta_{x, \varepsilon }$ such that $l := \ell
(\theta) > |\Gamma_1|$. Then there exist $i, j \in [1, l]$, $i <
j$, such that $\alpha_i (\theta) = \alpha_j (\theta)$. An easy
induction shows that $\alpha_{j + 1 - i} (\theta) = \alpha_1
(\theta) = \alpha_{x, \varepsilon}'$. Consequently, $\alpha_{j -
i} (\theta) \in \Sigma_{x, \varepsilon}$, which is impossible.
\end{proof}

Let $\omega$, $\sigma$, and $\omega'$ be as in the definition of
${}_+ \omega$. Obviously, $\alpha_{t \sigma, - T \sigma} =
\varnothing$, hence the above lemma implies that $\theta_{t
\sigma, - T \sigma} \neq \varnothing$. Now assume that $\ell
(\sigma) = 0$. This assumption means that $\alpha_{t \omega,
\varepsilon} = \varnothing$, where
\[
\varepsilon =
\begin{cases}
T \omega & \text{$r (\omega) > 0$},
\\
- T \omega & \text{$r (\omega) = 0$}.
\end{cases}
\]
Consequently, if $r (\omega) = 0$, i.e.\ $\omega' = \omega$, then
$\theta_{t \omega', - T \omega'} \neq \varnothing$. On the other
hand, if $r (\omega) > 0$, then $\theta_{t \omega, T \omega} \neq
\varnothing$. Moreover, in this case $\alpha_1 (\omega) \cdots
\alpha_{r (\omega)} (\omega) \theta \in \Theta_{t \omega, T
\omega}$ for each $\theta \in \Theta_{t \omega', - T \omega'}$,
hence we also get that $\theta_{t \omega', - T \omega'} \neq
\varnothing$.

The following lemma is crucial.

\begin{lemma} \label{lemma_AR}
Let $\omega$ be a homotopy band in $\bGamma$. Then ${}_+ \omega =
\varnothing$ if and only if ${}_+ (\psi \omega) = \varnothing$.
Moreover, if ${}_+ \omega \neq \varnothing$, then ${}_+ (\psi
\omega) = \Delta^{- m' (\omega)} (\psi ({}_+ \omega))$.
\end{lemma}

\begin{proof}
We make some remarks about the proof and leave the details to the
reader. Let $\sigma$ and $\omega'$ be as in the definition of
${}_+ \omega$.

First, if $\ell (\sigma) > 0$, then we prove $\Delta^{-1} ({}_+
(\psi \omega)) = \Delta^{- \ell (\theta_{t \sigma, - T \sigma})}
(\psi ({}_+ \omega))$. In the proof we use the following fact,
which can be easily proved by induction.

\begin{fact}
Let $\omega_0$ be a homotopy string in $\bGamma$ such that
$\sigma_{t \partial' (\psi \omega_0), - T \partial' (\psi
\omega_0)} = \sigma_{\omega_0} [0]$. Let $l \in [0, \ell
(\omega_0)]$. If $\sigma_i (\omega_0) \in \Gamma_1^{-1}$ for each
$i \in [1, l]$, then
\[
\Delta^{-l} (\psi \omega_0) = \sigma_{t \partial' (\psi
\omega_0'), - T \partial' (\psi \omega_0')} \cdot
\partial' (\psi \omega_0').
\]
where $\omega_0' := {}^{[l]} \omega$.
\end{fact}

Next, we assume that $\ell (\sigma) = 0$ and either $\ell
(\theta_{t \omega', - T \omega'}) > 0$ or $\ell (\omega') > 0$. In
this case we show either ${}_+ (\Delta^{r (\omega) - 1} (\psi
\omega)) = \Delta^{- \ell (\theta_{t \omega', - T \omega'})} (\psi
({}_+ \omega))$ if $\omega' \neq {}^{[r (\omega)]} \omega$, or
${}_+ (\Delta^{r (\omega)} (\psi \omega)) = \Delta^{- \ell
(\theta_{t \omega', - T \omega'}) + 1} (\psi ({}_+ \omega))$,
otherwise. In addition to the above fact, we use here also the
following.

\begin{fact}
Let $\omega_0$ be a homotopy string in $\bGamma$ such that $\ell
(\sigma_{\omega_0}) = 0$. Let $r \in [0, \ell (\omega_0)]$. If
$\sigma_i (\omega_0) \in \Gamma_1$ for each $i \in [1, r]$, then
\[
\Delta^r (\psi \omega_0) = \partial' (\psi \omega_0'),
\]
where $\omega_0' := {}^{[r]} \omega_0$.
\end{fact}

Finally, we prove that ${}_+ (\Delta^{r (\omega)} (\psi \omega)) =
\varnothing$ if $\ell (\sigma) = 0$, $\ell (\theta_{t \omega', - T
\omega'}') = 0$, and $\ell (\omega') = 0$. In this proof we also
use the latter fact.
\end{proof}

Dually, we put $\omega_+ := ({}_+ (\omega^{-1}))^{-1}$ for a
homotopy band $\omega$ in $\bGamma$. Lemma~\ref{lemma_AR} and
Corollary~\ref{coro_phiomegainverse} imply that $(\psi \omega)_+ =
\psi (\omega_+)$ for each homotopy band $\omega$ in $\bGamma$.
Finally, we put
\begin{gather*}
{}_+ \omega_+ :=
\begin{cases}
({}_+ \omega)_+ & \text{${}_+ \omega \neq \varnothing$},
\\
{}_+ (\omega_+) & \text{$\omega_+ \neq \varnothing$},
\end{cases}
\qquad \text{and} \qquad  m'' (\omega) :=
\begin{cases}
m' (\omega) & \text{${}_+ \omega \neq \varnothing$},
\\
m' (\omega_+) & \text{$\omega_+ \neq \varnothing$}.
\end{cases}
\end{gather*}

We obtain the following description of the almost split triangles
in $\calK^b (\bGamma)$ involving the string complexes, where
$X_{m, \varnothing}$ is the zero complex for $m \in \bbZ$.

\begin{maintheorem}[Part II: String complexes]
Let $\omega$ be a homotopy string in $\bGamma$ and $m \in \bbZ$.
Then we have an almost split triangle in $\calK^b (\bGamma)$ of
the form
\[
X_{m, \omega} \to X_{m + m' (\omega), {}_+ \omega} \oplus X_{m,
\omega_+} \to X_{m + m'' (\omega), {}_+ \omega_+} \to X_{m - 1,
\omega}.
\]
\end{maintheorem}

As a consequence we obtain the following description of the almost
split sequences with indecomposable middle terms containing the
string complexes (we encourage the reader to compare this result
with ~\cite{AvellaAlaminosGeiss2008}).

\begin{corollary}
Let $\omega$ be a homotopy string in $\bGamma$. Then $\varnothing
\in \{ {}_+ \omega, \omega_+ \}$ if and only if $\omega =
\theta_{x, \varepsilon}^{\varepsilon'}$ for $x \in \Gamma_0$ and
$\varepsilon, \varepsilon' \in \{ \pm 1 \}$ such that $\alpha_{x,
\varepsilon} = \varnothing$. Moreover, if this is the case, then
${}_+ \omega_+ = \theta_{t \sigma, - T \sigma}^{\varepsilon'}$,
where $\sigma := \sigma_{s \omega^{\varepsilon'}, - S
\omega^{\varepsilon'}}$.
\end{corollary}

\bibsection

\end{document}